\documentclass[a4paper,12pt,british,english]{article}
\usepackage[a4paper, margin=2cm]{geometry}
\usepackage[T1]{fontenc}
\usepackage[latin9]{inputenc}
\usepackage{color}
\usepackage{amsmath}
\usepackage{amssymb}
\usepackage{esint}
\PassOptionsToPackage{normalem}{ulem}
\usepackage{ulem}
\usepackage{enumerate}

\usepackage{babel}

\newtheorem{theo}{Theorem}

\newtheorem{prop}{Proposition}  
\newtheorem{coro}{Corollary}

\newtheorem{rem}{Remark}
\newtheorem{lema}{Lemma}
\newtheorem{fact}{Fact}
\newtheorem{defi}{Definition}

\newenvironment{dwd}{\par\noindent{\bf Proof.}}{\par\rightline{$\blacksquare$}}

 \global\long\def\sbr#1{\left[ #1\right] }
 \global\long\def\cbr#1{\left\{  #1\right\}  }
 \global\long\def\rbr#1{\left(#1\right)}

 \global\long\def\R{\mathbb{R}}

 \global\long\def\TDD#1{}
 \global\long\def\dd#1{\textnormal{d}#1}

 \global\long\def\TTV#1#2#3{\text{TV}^{#3}\!\rbr{#1,#2}}
 \global\long\def\thmTTV#1#2#3{\emph{TV}^{#3}\!\rbr{#1,#2}}

 \global\long\def\ra{\rightarrow}
 
 \global\long\def\ns{\infty}
 
  \global\long\def\TVnorm#1#2{\left\Vert {#1} \right\Vert_{{#2}-\text{TV},\left[a;b\right]}}
\global\long\def\TVnormthm#1#2{\left\Vert {#1} \right\Vert_{{#2}-\emph{TV},\left[a;b\right]}}

\global\long\def\Varnormthm#1#2{\left\Vert {#1} \right\Vert_{{#2}-\emph{var},\left[a;b\right]}}
\global\long\def\Oscnorm#1{\left\Vert {#1} \right\Vert_{\text{osc},\left[a;b\right]}}
\global\long\def\Oscnormthm#1{\left\Vert {#1} \right\Vert_{\emph{osc},\left[a;b\right]}}
\global\long\def\TVfullnorm#1#2{\left\Vert {#1} \right\Vert_{{TV},{#2},\left[a;b\right]}}
\global\long\def\TVfullnormA#1#2#3{\left\Vert {#1} \right\Vert_{TV,{#2},{#3}}}
\global\long\def\TVnormA#1#2#3{\left\Vert {#1} \right\Vert_{{#2}-\text{TV},{#3}}}
\global\long\def\TVnormAthm#1#2#3{\left\Vert {#1} \right\Vert_{{#2}-\emph{TV},{#3}}}

\begin{document}

\title{Integration of rough paths - the truncated variation approach}

\author{Rafa{\l{}} M. \L{}ochowski\\
E-mail : rlocho314@gmail.com}

\maketitle

\begin{abstract}
Using truncated variation techniques we obtain an improved version of the Lo\'{e}ve-Young inequality for the Riemann-Stieltjes integrals driven by rough paths. This allowed us to strenghten some result on the existence of solutions of integral equations driven by moderately irregular signals. 
We introduce also a new Banach space, containing as a proper subspace the paths with finite $p$-variation, and develop, in a systematic way, several parallel results for the paths from this space, obtained so far for the paths with finite $p$-variation. 
We start the paper with a general theorem on the existence of the Riemann-Stieltjes integral.
\end{abstract}

\section{Introduction}

The purpose of this paper is to investigate the top-down structure of the Riemann-Stieltjes integral and to state some general condition guaranteeing the existence of this integral, expressed in terms of the functional called {\em truncated variation}. 
We will also establish a quantitative relationship between the assymptotic behaviour of this functional (when its parameter called the {\em truncation parameter} tends to $0$) and the $p$-variation ($p\geq 1$). 

For $f:[a;b]\ra \R$  its truncated variation with the truncation parameter $\delta \geq 0$ will be denoted by $\TTV f{[a;b]}{\delta}.$ It may be simply defined as the greatest lower bound for the total variation of any function $g:[a;b]\ra \R,$ uniformly approximating $f$ with the accuracy $\delta/2,$
\[
\TTV f{[a;b]}{\delta} := \inf\cbr{\TTV {g}{[a;b]}{}: \|f - g \|_{\ns,[a;b]} \leq \delta/2}.
\]
$\|f - g \|_{\ns, [a;b]}$ denotes here $\sup_{a\leq t \leq b} \left|f(t) - g(t)\right|$ and the total variation $\TTV {g}{[a;b]}{}$ is defined as 
\[\TTV {g}{[a;b]}{} := \sup_{n}\sup_{a\leq t_{1}<t_{2}<\cdots<t_{n}\leq b}\sum_{i=1}^{n-1}\left|g\left(t_{i}\right)-g\left(t_{i-1}\right)\right|.  \]
It appears that the truncated variation $\TTV f{[a;b]}{\delta}$ is finite for any $\delta>0$ iff $f$ is {\em regulated} (cf. \cite[Fact 2.2]{LochowskiColloquium:2013}) and then for any $\delta >0$  the following equality holds
\begin{equation} \label{TV_def1}
\TTV f{[a;b]}{\delta} = \sup_{n}\sup_{a\leq t_{1}<t_{2}<\cdots<t_{n}\leq b}\sum_{i=1}^{n-1}\max\left\{ \left|f\left(t_{i}\right)-f\left(t_{i-1}\right)\right|-\delta,0\right\}
\end{equation}
(cf. \cite[Theorem 4]{LochowskiGhomrasniMMAS:2014}). The definition of the family of regulated functions is stated in the next section.
 
The simplest (and rather not interesting) case where the Riemann-Stieltjes integral  $\int_{a}^{b}f\left(t\right)\mathrm{d}g\left(t\right)$ (RSI in short) exists, is the situation when the integrand and integrator have no common points of discontinuity, the former is bounded and the latter has finite total variation. We will prove a general theorem (Theorem \ref{main}) encompassing this situation as well as a more interesting case when the integrand
and integrator have possibly unbounded variation, but they have finite
$p-$variation and $q-$variation respectively, with $p>1,$ $q>1$
and $p^{-1}+q^{-1}>1.$ The latter result is due to Young (\cite[Theorem on Stieltjes integrability, p. 264]{Young:1936}). For $f:\left[a;b\right]\rightarrow\R$
and $p>0,$ the $p-$variation, which we will denote by $V^{p}\left(f;\left[a;b\right]\right),$
is defined as 
\[
V^{p}\left(f,\left[a;b\right]\right)=\sup_{n}\sup_{a\leq t_{1}<t_{2}<\cdots<t_{n}\leq b}\sum_{i=1}^{n-1}\left|f\left(t_{i}\right)-f\left(t_{i-1}\right)\right|^{p}.
\]
The original Young's proof utlilizes elementary but clever induction
argument for finite sequences. Other proof of the Young theorem may
be found in \cite[Chapt. 6]{Friz:2010fk}, where integral estimates
based on control function and the Young-Lóeve inequality are used. These
approach is further applied in the rough-path theory setting. Further generalisations of Young's theorem are possible, with $p$-variation replaced by more general $\varphi$-variation: 
\[
V^{\varphi}\left(f,\left[a;b\right]\right)=\sup_{n}\sup_{a\leq t_{1}<t_{2}<\cdots<t_{n}\leq b}\sum_{i=1}^{n-1} \varphi \rbr{\left|f\left(t_{i}\right)-f\left(t_{i-1}\right)\right|},
\]
where $\varphi:[0;+\ns) \ra [0;+\ns)$ is a Young function, i.e. convex, strictly increasing function starting from $0$ (see e.g. \cite{Young:1938}, \cite{Dyackov:1988} and for a survey about another results of this type see the recent books \cite[Chapt. 3]{NorvaisaConcrete:2010},  \cite[Sect. 4.4]{Banas_et_al:2013}).

However, as far as we know, Theorem \ref{main} is a new result on the existence of the RSI. The proof of Theorem \ref{main} utilizes simple properties of the truncated variation and works for both (boundedness of integrand and finite total variation of integrator or finite $p$- and $q$- variations of integrand and integrator respectively) cases. This follows from the fact that the truncated variation with the truncation parameter $\delta = 0$ is simply the total variation. On the other hand there exists a simple relationship betwen the assymptotic behaviour of the truncated variation as $\delta \ra 0+$ and the finiteness of the $p$-variation.  
As a ``byproduct'' we will also obtain an inequality from which will follow a stronger version of the Lo\'{e}ve-Young inequality (see Corollaries \ref{corol_Young}, \ref{corol_Young1} and Remarks \ref{rem_Young}, \ref{rem_Young1}). As far as we know no variation of this inequality as in Corollaries \ref{corol_Young}, \ref{corol_Young1} and Remarks \ref{rem_Young}, \ref{rem_Young1} has yet appeared (see detailed historical notes on the the Lo\'{e}ve-Young inequality in \cite[pp. 212-214]{NorvaisaConcrete:2010}). We conjecture that using Theorem \ref{main} one may also obtain a variation of the  Lo\'{e}ve-Young inequality for $\varphi$-variation (see \cite[Theorem 3.89, Corollary 3.90]{NorvaisaConcrete:2010} or \cite[Theorem 4.40]{Banas_et_al:2013}). We intend to deal with this conjecture in the future. 

The already mentioned relationship betwen the assymptotic behaviour of the truncated variation as $\delta \ra 0+$ and the finiteness of the $p$-variation is not a completely new observation. Its qualitative version appears in an interesting paper by Tronel and Vladimirov, see \cite[Theorem 17]{TronelVladimirov:2000}. However, with the help of formula (\ref{TV_def1}) we will be able to obtain more precise, quantitative results. We will e.g. prove that for any $p\geq 1$ the following functional 
\begin{equation} \label{eq:TVnorm}
\| f\|_{TV,p,[a;b]}:= |f(a)|+ \sup_{\delta >0 } \rbr{ {\delta}^{p-1} \TTV f{[a;b]}{\delta} }^{1/p}
\end{equation}
defines a norm on the subspace of regulated functions $f:[a;b]\ra \R$ for which $\| f\|_{TV,p,[a;b]} < +\ns.$ Denote this space by ${\cal U}^p([a;b]).$ ${\cal U}^p([a;b])$ with norm (\ref{eq:TVnorm}) is a Banach space and contains, as a proper subspace, the space of functions with finite $p$-variation, ${\cal V}^p([a;b]).$ 
It is worth to mention that for example, for a typical path of a standard Brownian motion $B$ on $[0;T],$ $T>0,$ one has $V^2\left(B,\left[0;T\right]\right) = +\ns$ almost surely (see \cite{Levy:1940}), hence $B\notin {\cal V}^2([0;T])$ but $\| B \|_{TV,2,[0;T]} < +\ns$ almost surely, thus $B\in {\cal U}^2([0;T]) \supset {\cal V}^2([0;T]).$

Using Theorem \ref{main} we will also be able to state for the RSI $\int_a^b f(t) \dd g(t)$ an inequality of the Lo\'{e}ve-Young type, whenever $f \in {\cal U}^p([a;b])$ and $g \in {\cal U}^q([a;b])$ with $p>1,$ $q>1$ and $p^{-1}+q^{-1}>1.$
Next, following Lyons \cite{Lyons:1994}, we will establish an inequality for the $\| \cdot\|_{TV,p,[a;b]}$ norm of the function $[a;b] \ni t \mapsto\int_a^{t} f(s) \dd g(s)$ in terms of the norms $\| f\|_{TV,p,[a;b]}$ and $\| g\|_{TV,q,[a;b]}.$ The method of the proof will be completely different from \cite{Lyons:1994}, since we will not be able to use the superadditivity property which holds for the $p$-variation: for any $d\in(a;b),$ $V^p\left(f,\left[a;b\right]\right)  \geq  V^p\left(f,\left[a;d\right]\right) + V^p\left(f,\left[d;b\right]\right).$

After having obtained for the paths from the spaces ${\cal U}^p([a;b]),$ $p\geq 1,$ parallel results to the results obtained so far for the paths with finite $p$-variation, we will be able, following Lyons \cite{Lyons:1994}, \cite{LyonsCaruana:2007}, to solve few types of integral equations driven by moderately irregular signals from these spaces. By moderately irregular signals we mean continuous signals belonging to ${\cal U}^p([a;b])$ with $p\in (1;2).$ It is well known that for higher degrees of irregularity, corresponding to $p\geq 2,$ one needs, constructing approximations of integral equations, to consider terms of a new type (like L\'{e}vy's area). We believe that the tuncated variation approach for such paths is also possible and this will be a topic of our further research. 

Let us comment shortly on the organisation of the paper. In the next section we prove a general theorem on the existence of the Riemann-Stieltjes integral, expressed in terms of the truncated variation functionals and derive from it the Young theorem as well as the Lo\'{e}ve-Young inequality. Next, in Section 3, we prove that the space ${\cal U}^p([a;b]),$ $p\geq1,$ equipped with the norm (\ref{eq:TVnorm}) is a Banach space and then prove the Lyons type estimate for the norms of the integrals driven by signals from this space. In the last, fourth section we deal with  few types of integral equations driven by moderately irregular signals from the space ${\cal U}^p([a;b]).$
 
\section{A theorem on the existence of the Riemann-Stieltjes integral} \label{RSI_existence}

In this section we will prove a general theorem on the existence of the Riemann-Stieltjes integral $\int_{a}^{b}f\left(t\right)\mathrm{d}g\left(t\right)$ formulated in terms of the truncated variation.
We will assume that both - integrand $f:[a;b]\ra \R$ and integrator $g:[a;b]\ra \R$ are regulated functions. Let us recall the definition of a regulated function. 
\begin{defi}
A function  $h:[a;b]\ra \R$ is {\em regulated} if there exist one-sided limits
$\lim_{t\rightarrow a+}h\left(t\right)$ and $\lim_{t\rightarrow b-}h\left(t\right),$
and for any $t\in\left(a;b\right)$ there exist one-sided limits $\lim_{t\rightarrow x-}h\left(t\right)$
and $\lim_{t\rightarrow x+}h\left(t\right).$
\end{defi}
We will also need the following
result (cf. \cite[Theorem 4]{LochowskiGhomrasniMMAS:2014}): for any regulated
function $f:\left[a;b\right]\rightarrow\R$ and $\delta>0$ there
exists a regulated function $f^{\delta}:\left[a;b\right]\rightarrow\R$
such that $\left\Vert f-f^{\delta}\right\Vert _{\infty,\left[a;b\right]}\leq\delta/2$
and 
\[
\TTV{f^{\delta}}{\left[a;b\right]}0=\TTV f{\left[a;b\right]}{\delta}.
\]
Directly from the definition it follows that the truncated variation
is a superadditive functional of the interval, i.e. for any $d\in\left(a;b\right)$
\[
\TTV f{\left[a;b\right]}{\delta}\geq\TTV f{\left[a;d\right]}{\delta}+\TTV f{\left[d;b\right]}{\delta}.
\]
Moreover, we also have the following easy estimate of the truncated
variation of a function $f$ perturbed by some other function $h:$
\begin{equation}
\TTV{f+h}{\left[a;b\right]}{\delta}\leq\TTV f{\left[a;b\right]}{\delta}+\TTV h{\left[a;b\right]}0,\label{eq:TV_variation}
\end{equation}
which stems directly from the inequality: for $a\leq s<t\leq b,$
\begin{align*}
 & \max\left\{ \left|f\left(t\right)+h\left(t\right)-\left\{ f\left(s\right)+h\left(s\right)\right\} \right|-\delta,0\right\} \\
 & \leq\max\left\{ \left|f\left(t\right)-f\left(s\right)\right|-\delta,0\right\} +\left|h\left(t\right)-h\left(s\right)\right|.
\end{align*}

\begin{theo} \label{main} Let $f,g:\left[a;b\right]\rightarrow\R$
be two regulated functions which have no common points of discontinuity.
Let $\eta_{0}\geq\eta_{1}\geq\ldots$ and $\theta_{0}\geq\theta_{1}\geq\ldots$
be two sequences of non-negative numbers, such that $\eta_{k}\downarrow0,$
$\theta_{k}\downarrow0$ as $k\rightarrow+\infty.$ Define $\eta_{-1}:=\sup_{a\leq t\leq b}\left|f\left(t\right)-f\left(a\right)\right|$
and 
\[
S:=\sum_{k=0}^{+\infty}2^{k}\eta_{k-1}\cdot\TTV g{\left[a;b\right]}{\theta_{k}}+\sum_{k=0}^{\infty}2^{k}\theta_{k}\cdot\TTV f{\left[a;b\right]}{\eta_{k}}.
\]
If $S<+\infty$ then the Riemann-Stieltjes integral 
$
\int_{a}^{b}f\left(t\right)\mathrm{d}g\left(t\right)
$
exists and one has the following estimate 
\begin{equation}
\left|\int_{a}^{b}f\left(t\right)\mathrm{d}g\left(t\right)-f\left(a\right)\left[g\left(b\right)-g\left(a\right)\right]\right|\leq S.\label{eq:estimate_integral}
\end{equation}
\end{theo}

\begin{rem} The assumption that $f$ and $g$ has no common points
of discontinuity is necessary for the existence of the Riemann-Stieltjes
integral $\int_{a}^{b}f\left(t\right)\mathrm{d}g\left(t\right).$
When a more general integrals are considered (e.g. the Moore-Pollard
integral, c.f. \cite[p. 263]{Young:1936}), we may weaken this assumption
and assume that $f$ and $g$ have no common one-sided discontinuities.
\end{rem} The proof of Theorem \ref{main} will be based on the following
lemma.

\begin{lema} \label{lema} Let $f,g:\left[a;b\right]\rightarrow\R$
be two regulated functions. Let $c=t_{0}<t_{1}<\ldots<t_{n}=d$ be
any partition of the interval $\left[c;d\right]\subset\left[a;b\right]$
and let $\xi_{0}=c$ and $\xi_{1},\ldots,\xi_{n}$ be such that $t_{i-1}\leq\xi_{i}\leq t_{i}$
for $i=1,2,\ldots,n.$ Then for $\delta_{-1}:=\sup_{c\leq t\leq d}\left|f\left(t\right)-f\left(c\right)\right|,$
 $\delta_{0}\geq\delta_{1}\geq\ldots\geq\delta_{r}\geq0$ and
$\varepsilon_{0}\geq\varepsilon_{1}\geq\ldots\geq\varepsilon_{r}\geq0$
the following estimate holds 
\begin{eqnarray*}
 &  & \left|\sum_{i=1}^{n}f\left(\xi_{i}\right)\left[g\left(t_{i}\right)-g\left(t_{i-1}\right)\right]-f\left(c\right)\left[g\left(d\right)-g\left(c\right)\right]\right|\\
 &  & \leq\sum_{k=0}^{r}2^{k}\delta_{k-1}\cdot\TTV g{\left[c;d\right]}{\varepsilon_{k}}+\sum_{k=0}^{r}2^{k}\varepsilon_{k}\cdot\TTV f{\left[c;d\right]}{\delta_{k}}+n\delta_{r}\varepsilon_{r}.
\end{eqnarray*}
\end{lema} 
\begin{dwd} Denote $\varepsilon=\varepsilon_{0},$ by
summation by parts, we have the following equality 
\begin{align}
 & \sum_{i=1}^{n}  f\left(\xi_{i}\right)\left[g\left(t_{i}\right)-g\left(t_{i-1}\right)\right]-f\left(c\right)\left[g\left(d\right)-g\left(c\right)\right]  \nonumber \\
 & =  \sum_{i=1}^{n}\left[f\left(\xi_{i}\right)-f\left(c\right)\right]\left[g^{\varepsilon}\left(t_{i}\right)-g^{\varepsilon}\left(t_{i-1}\right)\right]\nonumber \\
 &  \quad+\sum_{i=1}^{n}\left[f\left(\xi_{i}\right)-f\left(c\right)\right]\left[g\left(t_{i}\right)-g^{\varepsilon}\left(t_{i}\right)-\left\{ g\left(t_{i-1}\right)-g^{\varepsilon}\left(t_{i-1}\right)\right\} \right]\nonumber \\
 & =  \sum_{i=1}^{n}\left[f\left(\xi_{i}\right)-f\left(c\right)\right]\left[g^{\varepsilon}\left(t_{i}\right)-g^{\varepsilon}\left(t_{i-1}\right)\right]\nonumber \\
 & \quad +\sum_{i=1}^{n}\left[g\left(d\right)-g^{\varepsilon}\left(d\right)-\left\{ g\left(t_{i-1}\right)-g^{\varepsilon}\left(t_{i-1}\right)\right\} \right]\left[f\left(\xi_{i}\right)-f\left(\xi_{i-1}\right)\right],\label{eq:jeden}
\end{align}
where $g^{\varepsilon}:[c;d]\ra \R$ is regulated and such that 
\[
\left\Vert g-g^{\varepsilon}\right\Vert _{\infty,[c;d]}\leq\frac{1}{2}\varepsilon\mbox{ and }\TTV{g^{\varepsilon}}{\left[c;d\right]}0=\TTV g{\left[c;d\right]}{\varepsilon}.
\]
Similarly, for $\delta=\delta_{0}$ we may write 
\begin{align}
 & \sum_{i=1}^{n} \left[g\left(d\right)-g^{\varepsilon}\left(d\right)-\left\{ g\left(t_{i-1}\right)-g^{\varepsilon}\left(t_{i-1}\right)\right\} \right]\left[f\left(\xi_{i}\right)-f\left(\xi_{i-1}\right)\right]\nonumber \\
 & =  \sum_{i=1}^{n}\left[g\left(d\right)-g^{\varepsilon}\left(d\right)-\left\{ g\left(t_{i-1}\right)-g^{\varepsilon}\left(t_{i-1}\right)\right\} \right]\left[f^{\delta}\left(\xi_{i}\right)-f^{\delta}\left(\xi_{i-1}\right)\right]\label{eq:dwa}\\
 &  \quad +\sum_{i=1}^{n}\left[f\left(\xi_{i}\right)-f^{\delta}\left(\xi_{i}\right)-\left\{ f\left(c\right)-f^{\delta}\left(c\right)\right\} \right]\left[\left\{ g\left(t_{i}\right)-g^{\varepsilon}\left(t_{i}\right)\right\} -\left\{ g\left(t_{i-1}\right)-g^{\varepsilon}\left(t_{i-1}\right)\right\} \right],\nonumber 
\end{align}
where $f^{\delta}:[c;d]\ra \R$ is regulated and such that 
\[
\left\Vert f-f^{\delta}\right\Vert _{\infty,[c;d]}\leq\frac{1}{2}\delta\mbox{ and }\TTV{f^{\delta}}{\left[c;d\right]}0=\TTV f{\left[c;d\right]}{\delta}.
\]
Since $\TTV{g^{\varepsilon}}{\left[c;d\right]}0=\TTV g{\left[c;d\right]}{\varepsilon},$
$\TTV{f^{\delta}}{\left[c;b\right]}0=\TTV f{\left[c;d\right]}{\delta},$
$\left\Vert g-g^{\varepsilon}\right\Vert _{\infty,[c;d]}\leq\varepsilon/2$
and $\left\Vert f-f^{\delta}\right\Vert _{\infty,[c;d]}\leq\delta/2,$
from (\ref{eq:jeden}) and (\ref{eq:dwa}) we have the following estimate
\begin{eqnarray}
 &  & \left|\sum_{i=1}^{n}f\left(\xi_{i}\right)\left[g\left(t_{i}\right)-g\left(t_{i-1}\right)\right]-f\left(c\right)\left[g\left(d\right)-g\left(c\right)\right]\right|\nonumber \\
 &  & \leq\sup_{c\leq t\leq d}\left|f\left(t\right)-f\left(c\right)\right|\cdot\TTV g{\left[c;d\right]}{\varepsilon}+\varepsilon\cdot\TTV f{\left[c;d\right]}{\delta}+n\delta\varepsilon.\label{eq:jeden1}
\end{eqnarray}
Denote $g_{1}:=g-g^{\varepsilon},$ $f_{1}:=f-f^{\delta}$ on $[c;d].$ By (\ref{eq:jeden})
and (\ref{eq:dwa}), instead of the last summand $n\delta\varepsilon$
in (\ref{eq:jeden1}) we may write the estimate 
\begin{eqnarray}
 &  & \left|\sum_{i=1}^{n}\left[f\left(\xi_{i}\right)-f^{\delta}\left(\xi_{i}\right)-\left\{ f\left(c\right)-f^{\delta}\left(c\right)\right\} \right]\left[\left\{ g\left(t_{i}\right)-g^{\varepsilon}\left(t_{i}\right)\right\} -\left\{ g\left(t_{i-1}\right)-g^{\varepsilon}\left(t_{i-1}\right)\right\} \right]\right|\nonumber \\
 &  & =\left|\sum_{i=1}^{n}\left[f_{1}\left(\xi_{i}\right)-f_{1}\left(c\right)\right]\left[g_{1}\left(t_{i}\right)-g_{1}\left(t_{i-1}\right)\right]\right|\nonumber \\
 &  & \leq\sup_{a\leq t\leq b}\left|f_{1}\left(t\right)-f_{1}\left(c\right)\right|\cdot\TTV{g_{1}}{\left[c;d\right]}{\varepsilon_{1}}+\varepsilon_{1}\cdot\TTV{f_{1}}{\left[c;d\right]}{\delta_{1}}+n\delta_{1}\varepsilon_{1}\nonumber \\
 &  & \leq\delta\cdot\TTV{g_{1}}{\left[c;d\right]}{\varepsilon_{1}}+\varepsilon_{1}\cdot\TTV{f_{1}}{\left[c;d\right]}{\delta_{1}}+n\delta_{1}\varepsilon_{1},\label{eq:dwa1}
\end{eqnarray}
where the last but one inequlity in (\ref{eq:dwa1}) follows by the
same reasoning for $f_{1}$ and $g_{1}$ as inequality (\ref{eq:jeden})
for $f$ and $g.$ Repeating these arguments, by induction we get
\begin{eqnarray}
 &  & \left|\sum_{i=1}^{n}f\left(\xi_{i}\right)\left[g\left(t_{i}\right)-g\left(t_{i-1}\right)\right]-f\left(c\right)\left[g\left(d\right)-g\left(c\right)\right]\right|\nonumber \\
 &  & \leq\sum_{k=0}^{r}\delta_{k-1}\cdot\TTV{g_{k}}{\left[c;d\right]}{\varepsilon_{k}}+\sum_{k=0}^{r}\varepsilon_{k}\cdot\TTV{f_{k}}{\left[c;d\right]}{\delta_{k}}+n\delta_{r}\varepsilon_{r},\label{eq:induction}
\end{eqnarray}
where $\delta_{-1}:=\sup_{c\leq t\leq c}\left|f\left(t\right)-f\left(a\right)\right|,$
$g_{0}\equiv g,$ $f_{0}\equiv f$ and for $k=1,2,\ldots,r,$ $g_{k}:=g_{k-1}-g_{k-1}^{\varepsilon_{k-1}},$
$f_{k}:=f_{k-1}-f_{k-1}^{\delta_{k-1}}$ are defined similarly as
$g_{1}$ and $f_{1}.$\par Since $\varepsilon_{k}\leq\varepsilon_{k-1}$
for $k=1,2,\ldots,r,$ by (\ref{eq:TV_variation}) and the fact that
the function $\delta\mapsto\TTV h{\left[c;d\right]}{\delta}$ is non-increasing,
we estimate 
\begin{eqnarray*}
\TTV{g_{k}}{\left[c;d\right]}{\varepsilon_{k}} & = & \TTV{g_{k-1}-g_{k-1}^{\varepsilon_{k-1}}}{\left[c;d\right]}{\varepsilon_{k}}\\
 & \leq & \TTV{g_{k-1}}{\left[c;d\right]}{\varepsilon_{k}}+\TTV{g_{k-1}^{\varepsilon_{k-1}}}{\left[c;d\right]}0\\
 & = & \TTV{g_{k-1}}{\left[c;d\right]}{\varepsilon_{k}}+\TTV{g_{k-1}}{\left[c;d\right]}{\varepsilon_{k-1}}\\
 & \leq & 2\TTV{g_{k-1}}{\left[c;d\right]}{\varepsilon_{k}}.
\end{eqnarray*}
Hence, by recursion, for $k=1,2,\ldots,r,$ 
\[
\TTV{g_{k}}{\left[c;d\right]}{\varepsilon_{k}}\leq2^{k}\TTV g{\left[c;d\right]}{\varepsilon_{k}}.
\]
Similarly, for $k=1,2,\ldots,r,$ we have 
\[
\TTV{f_{k}}{\left[c;d\right]}{\delta_{k}}\leq2^{k}\TTV f{\left[c;d\right]}{\delta_{k}}.
\]
By (\ref{eq:induction}) and last two estimates we get the desired
estimate.
\end{dwd} 
\begin{rem} \label{symmetry}
Notice that starting in (\ref{eq:jeden}) from the summation by parts, then splitting the difference $f\left(\xi_{i}\right)-f\left(\xi_{i-1}\right):$
\begin{align*}
 & \sum_{i=1}^{n}  f\left(\xi_{i}\right)\left[g\left(t_{i}\right)-g\left(t_{i-1}\right)\right]-f\left(c\right)\left[g\left(d\right)-g\left(c\right)\right]  \\
 & =  \sum_{i=1}^{n}\left[g\left(d\right)-g\left(t_{i-1}\right)\right]
\left[f\left(\xi_{i}\right)-f\left(\xi_{i-1}\right)\right]
\\
 & =  \sum_{i=1}^{n}\left[g\left(d\right)-g\left(t_{i-1}\right)\right]\left[f^{\delta}\left(\xi_{i}\right)-f^{\delta}\left(\xi_{i-1}\right)\right] \\
& \quad + \sum_{i=1}^{n}\left[g\left(d\right)-g\left(t_{i-1}\right)\right]\left[f\left(\xi_{i}\right)-f^{\delta}\left(\xi_{i}\right)
-\cbr{f\left(\xi_{i-1}\right) - f^{\delta}\left(\xi_{i-1}\right) }\right] 
\end{align*}
and proceeding similarly as in the proof of Lemma \ref{lema} we get the symmetric estimate 
\begin{eqnarray}
 &  & \left|\sum_{i=1}^{n}f\left(\xi_{i}\right)\left[g\left(t_{i}\right)-g\left(t_{i-1}\right)\right]-f\left(c\right)\left[g\left(d\right)-g\left(c\right)\right]\right| \nonumber \\
 &  & \leq\sum_{k=0}^{r}2^{k}\varepsilon_{k-1}\cdot\TTV f{\left[c;d\right]}{\delta_{k}}+\sum_{k=0}^{r}2^{k}\delta_{k}\cdot\TTV g{\left[c;d\right]}{\varepsilon_{k}}+n\delta_{r}\varepsilon_{r}, \label{symmetry_left}
\end{eqnarray}
where $\varepsilon_{-1}=\sup_{c\leq t \leq d}\left|g(d)-g(t) \right|.$
\end{rem}
\begin{rem} \label{n=1}
Setting in Lemma \ref{lema}, $n=1$ for any $\xi \in [c;d]$ we get the estimate 
\begin{eqnarray}
 &  & \left|\left(f\left(\xi\right)-f\left(c\right)\right)\left[g\left(d\right)-g\left(c\right)\right]\right| \nonumber \\
 &  & \leq\sum_{k=0}^{r}2^{k}\delta_{k-1}\cdot\TTV f{\left[c;d\right]}{\varepsilon_{k}}+\sum_{k=0}^{r}2^{k}\varepsilon_{k}\cdot\TTV g{\left[c;d\right]}{\delta_{k}}+n\delta_{r}\varepsilon_{r}. \label{n=1_left}
\end{eqnarray}
and similarly, setting in Remark \ref{symmetry}, $n=1$ we get similar estimate, where the left side of (\ref{n=1_left}) is replaced by the left side of  (\ref{symmetry_left}).
\end{rem}
Now we proceed to the proof of Theorem \ref{main}.
\begin{dwd} It is enough to prove that for any two partitions 
\[
\pi=\left\{ a=a_{0}<a_{1}<\ldots<a_{l}=b\right\} ,
\]
\[
\rho=\left\{ a=b_{0}<b_{1}<\ldots<b_{m}=b\right\} 
\]
and $\nu_{i}\in\left[a_{i-1};a_{i}\right],$ $\xi_{j}\in\left[b_{j-1};b_{j}\right],$
$i=1,2,\ldots,l,$ $j=1,2,\ldots,m,$ the difference 
\[
\left|\sum_{i=1}^{l}f\left(\nu_{i}\right)\left[g\left(a_{i}\right)-g\left(a_{i-1}\right)\right]-\sum_{j=1}^{m}f\left(\xi_{j}\right)\left[g\left(b_{j}\right)-g\left(b_{j-1}\right)\right]\right|
\]
is as small as we please, provided that the meshes of the partitions
$\pi$ and $\rho$, defined as 
\[
\mbox{mesh}\left(\pi\right):=\max_{i=1,2,\ldots,l}\left(a_{i}-a_{i-1}\right),\text{ }\mbox{mesh}\left(\rho\right):=\max_{j=1,2,\ldots,m}\left(b_{j}-b_{j-1}\right)
\]
respectively, are sufficiently small. Define 
\[
\sigma=\pi\cup\rho=\left\{ a=s_{0}<s_{1}<\ldots<s_{n}=b\right\} 
\]
and for $i=1,2,\ldots,l$ consider 
\[
\left|f\left(\nu_{i}\right)\left[g\left(a_{i}\right)-g\left(a_{i-1}\right)\right]-\sum_{k:s_{k-1},s_{k}\in\left[a_{i-1};a_{i}\right]}f\left(s_{k-1}\right)\left[g\left(s_{k}\right)-g\left(s_{k-1}\right)\right]\right|.
\]
We estimate 
\begin{eqnarray*}
\lefteqn{\left|f\left(\nu_{i}\right)\left[g\left(a_{i}\right)-g\left(a_{i-1}\right)\right]-\sum_{k:s_{k-1},s_{k}\in\left[a_{i-1};a_{i}\right]}f\left(s_{k-1}\right)\left[g\left(s_{k}\right)-g\left(s_{k-1}\right)\right]\right|}\\
 & \leq & \left|f\left(\nu_{i}\right)\left[g\left(a_{i}\right)-g\left(a_{i-1}\right)\right]-f\left(a_{i-1}\right)\left[g\left(a_{i}\right)-g\left(a_{i-1}\right)\right]\right|\\
 &  & +\left|\sum_{k:s_{k-1},s_{k}\in\left[a_{i-1};a_{i}\right]}f\left(s_{k-1}\right)\left[g\left(s_{k}\right)-g\left(s_{k-1}\right)\right]-f\left(a_{i-1}\right)\left[g\left(a_{i}\right)-g\left(a_{i-1}\right)\right]\right|.
\end{eqnarray*}

Recall the definition of  $S.$ If there exists $N=0,1,2,\ldots$ such that $\eta_{N}=0$ or $\theta_{N}=0$ then $\TTV f{\left[a;b\right]}{}$ or $\TTV g{\left[a;b\right]}{}$ is finite, moreover, both functions, $f$ and $g,$ are bounded (since they are regulated), hence the integral $\int_a^b f(t)\dd g(t)$ exists. Thus we may and will assume that 
$\eta_{N}>0$ and $\theta_{N}>0$ for all $N=0,1,2,\ldots$

Choose $N=1,2,\ldots.$ By the assumption that $f$ and $g$ have
no common points of discontinuity, for sufficiently small $\mbox{mesh}\left(\pi\right),$
for $i=1,2,\ldots,l$ we have 
\begin{equation}
\sup_{a_{i-1}\leq s\leq a_{i}}\left|f\left(s\right)-f\left(a_{i-1}\right)\right|\leq\eta_{N-1}\label{eq:ind_f}
\end{equation}
or 
\begin{equation}
\sup_{a_{i-1}\leq s\leq a_{i}}\left|g\left(a_{i}\right)-g\left(s\right)\right|\leq\theta_{N-1}.\label{eq:ind_J}
\end{equation}
To see this, assume that for every $h>0,$ there exist $\left[a_{h};b_{h}\right]\subset\left[a;b\right]$
such that $b_{h}-a_{h}\leq h$ and $\sup_{x,y\in\left[a_{h};b_{h}\right]}\left|f\left(y\right)-f\left(y\right)\right|>\eta_{N-1}$
and $\sup_{x,y\in\left[a_{h};b_{h}\right]}\left|g\left(x\right)-g\left(y\right)\right|>\theta_{N-1}.$
We choose a convergent subsequence of the sequence $\left(a_{1/n}+b_{1/n}\right)/2,$
$n=1,2,.\ldots,$ and we see that the limit of this sequence is a
point of discontinuity for both $f$ and $g,$ which is a contradiction
with the assumption that $f$ and $g$ have no common points of discontinuity.\par Let
$I$ be the set of all indices $i=1,2,\ldots,l$ for which (\ref{eq:ind_f})
holds. Now, for $i\in I,$ set $\delta_{j-1}:=\eta_{N+j-1},$ $\varepsilon_{j}:=\theta_{N+j},$ $j=0,1,2,\ldots,$ and 
define 
\[
S_{i}:=\sum_{j=0}^{+\infty}2^{j}\eta_{j-1}\cdot\TTV g{\left[a_{i-1};a_{i}\right]}{\theta_{j}}+\sum_{j=0}^{+\infty}2^{j}\theta_{j}\cdot\TTV f{\left[a_{i-1};a_{i}\right]}{\eta_{j}}.
\]
By Lemma \ref{lema} we estimate 
\begin{eqnarray*}
 &  & \left|\sum_{k:s_{k-1},s_{k}\in\left[a_{i-1};a_{i}\right]}f\left(s_{k-1}\right)\left[g\left(s_{k}\right)-g\left(s_{k-1}\right)\right]-f\left(a_{i-1}\right)\left[g\left(a_{i}\right)-g\left(a_{i-1}\right)\right]\right|\\
 &  & \leq\sum_{j=0}^{+\infty}2^{j}\delta_{j-1}\cdot\TTV g{\left[a_{i-1};a_{i}\right]}{\varepsilon_{j}}+\sum_{j=0}^{+\infty}2^{j}\varepsilon_{j}\cdot\TTV f{\left[a_{i-1};a_{i}\right]}{\delta_{j}}\\
 &  & \leq\sum_{j=0}^{+\infty}2^{j}\eta_{N+j-1}\cdot\TTV g{\left[a_{i-1};a_{i}\right]}{\theta_{N+j}}+\sum_{j=0}^{+\infty}2^{j}\theta_{N+j}\cdot\TTV f{\left[a_{i-1};a_{i}\right]}{\eta_{N+j}}\\
 &  & \leq2^{-N}S_{i}.
\end{eqnarray*}
Similarly, 
\begin{eqnarray*}
\left|f\left(\nu_{i}\right)\left[g\left(a_{i}\right)-g\left(a_{i-1}\right)\right]-f\left(a_{i-1}\right)\left[g\left(a_{i}\right)-g\left(a_{i-1}\right)\right]\right| & \leq & 2^{-N}S_{i}.
\end{eqnarray*}
Hence 
\begin{eqnarray}
 &  & \left|f\left(\nu_{i}\right)\left[g\left(a_{i}\right)-g\left(a_{i-1}\right)\right]-\sum_{k:s_{k-1},s_{k}\in\left[a_{i-1};a_{i}\right]}f\left(s_{k-1}\right)\left[g\left(s_{k}\right)-g\left(s_{k-1}\right)\right]\right|\nonumber \\
 &  & \leq2^{1-N}S_{i}.\label{eq:ineq_i}
\end{eqnarray}
The truncated variation is a superadditive functional of the interval, from
which we have 
\[
\sum_{i\in I}\TTV g{\left[a_{i-1};a_{i}\right]}{\theta_{j}}\leq\TTV g{\left[a;b\right]}{\theta_{j}},
\]
\[
\sum_{i\in I}\TTV f{\left[a_{i-1};a_{i}\right]}{\eta_{j}}\leq\TTV f{\left[a;b\right]}{\eta_{j}}.
\]
By (\ref{eq:ineq_i}) and last two inequalities, summing over $i\in I$
we get the estimate 
\begin{eqnarray}
 &  & \left|\sum_{i\in I}\left\{ f\left(\nu_{i}\right)\left[g\left(a_{i}\right)-g\left(a_{i-1}\right)\right]-\sum_{k:s_{k-1},s_{k}\in\left[a_{i-1};a_{i}\right]}f\left(s_{k-1}\right)\left[g\left(s_{k}\right)-g\left(s_{k-1}\right)\right]\right\} \right|\nonumber \\
 &  & \leq2^{1-N}\sum_{i\in I}S_{i}\leq2^{1-N}S.\label{eq:estim_I-1}
\end{eqnarray}
\par Now, let $J$ be the set of all indices, for which (\ref{eq:ind_J})
holds. For $i=1,2,\ldots,l$ define 
\[
T_{i}:=\sum_{j=0}^{+\infty}2^{j}\theta_{j}\cdot\TTV f{\left[a_{i-1};a\right]}{\eta_{j}}+\sum_{j=0}^{+\infty}2^{j}\eta_{j}\cdot\TTV g{\left[a_{i-1};a_{i}\right]}{\theta_{j+1}}.
\]
For $i\in J,$ by the summation by parts and then by Lemma \ref{lema}
we get\par 
\begin{eqnarray*}
 &  & \left|f\left(a_{i}\right)\left[g\left(a_{i}\right)-g\left(a_{i-1}\right)\right]-\sum_{k:s_{k-1},s_{k}\in\left[a_{i-1};a_{i}\right]}f\left(s_{k-1}\right)\left[g\left(s_{k}\right)-g\left(s_{k-1}\right)\right]\right|\\
 &  & =\left|\sum_{k:s_{k-1},s_{k}\in\left[a_{i-1};a_{i}\right]}g\left(s_{k}\right)\left[f\left(s_{k}\right)-f\left(s_{k-1}\right)\right]-g\left(a_{i-1}\right)\left[f\left(a_{i}\right)-f\left(a_{i-1}\right)\right]\right|\\
 &  & \leq\sum_{j=0}^{+\infty}2^{j}\theta_{N+j-1}\cdot\TTV f{\left[a_{i-1};a_{i}\right]}{\eta_{N+j}}+\sum_{j=0}^{+\infty}2^{j}\eta_{N+j}\cdot\TTV g{\left[a_{i-1};a_{i}\right]}{\theta_{N+j}}.\\
 &  & \leq2^{1-N}T_{i}\leq2^{1-N}S_{i}.
\end{eqnarray*}
Similarly, by Lemma \ref{lema}, 
\begin{eqnarray*}
 &  & \left|f\left(a_{i}\right)\left[g\left(a_{i}\right)-g\left(a_{i-1}\right)\right]-f\left(\nu_{i}\right)\left[g\left(a_{i}\right)-g\left(a_{i-1}\right)\right]\right|\\
 &  & =\left|g\left(a_{i-1}\right)\left[f\left(\nu_{i}\right)-f\left(a_{i-1}\right)\right]+g\left(a_{i}\right)\left[f\left(a_{i}\right)-f\left(\nu_{i}\right)\right]-g\left(a_{i-1}\right)\left[f\left(a_{i}\right)-f\left(a_{i-1}\right)\right]\right|\\
 &  & \leq2^{1-N}S_{i}.
\end{eqnarray*}
From last two inequalities we get 
\begin{eqnarray*}
\lefteqn{\left|f\left(\nu_{i}\right)\left[g\left(a_{i}\right)-g\left(a_{i-1}\right)\right]-\sum_{k:s_{k-1},s_{k}\in\left[a_{i-1};a_{i}\right]}f\left(s_{k-1}\right)\left[g\left(s_{k}\right)-g\left(s_{k-1}\right)\right]\right|}\\
 & \leq & \left|f\left(a_{i}\right)\left[g\left(a_{i}\right)-g\left(a_{i-1}\right)\right]-f\left(\nu_{i}\right)\left[g\left(a_{i}\right)-g\left(a_{i-1}\right)\right]\right|\\
 &  & +\left|f\left(a_{i}\right)\left[g\left(a_{i}\right)-g\left(a_{i-1}\right)\right]-\sum_{k:s_{k-1},s_{k}\in\left[a_{i-1};a_{i}\right]}f\left(s_{k-1}\right)\left[g\left(s_{k}\right)-g\left(s_{k-1}\right)\right]\right|\\
 & \leq & 2^{2-N}S_{i}.
\end{eqnarray*}
Summing over $i\in J$ and using the superadditivity of the truncated
variation as a function of the interval, we get the estimate 
\begin{eqnarray}
 &  & \left|\sum_{i\in J}\left\{ f\left(\nu_{i}\right)\left[g\left(a_{i}\right)-g\left(a_{i-1}\right)\right]-\sum_{k:s_{k-1},s_{k}\in\left[a_{i-1};a_{i}\right]}f\left(s_{k-1}\right)\left[g\left(s_{k}\right)-g\left(s_{k-1}\right)\right]\right\} \right|\nonumber \\
 &  & \leq2^{2-N}\sum_{i\in J}S_{i}\leq2^{2-N}S.\label{eq:estim_J-1}
\end{eqnarray}
Finally, from (\ref{eq:estim_I-1}) and (\ref{eq:estim_J-1}) we get
\[
\left|f\left(\nu_{i}\right)\left[g\left(b\right)-g\left(a\right)\right]-\sum_{k=1}^{n}f\left(s_{k-1}\right)\left[g\left(s_{k}\right)-g\left(s_{k-1}\right)\right]\right|\leq6\cdot2^{-N}S.
\]
\par Similar estimate holds for 
\[
\left|\sum_{i=j}^{m}f\left(\xi_{j}\right)\left[g\left(b_{j}\right)-g\left(b_{j-1}\right)\right]-\sum_{k=1}^{n}f\left(s_{k-1}\right)\left[g\left(s_{k}\right)-g\left(s_{k-1}\right)\right]\right|,
\]
provided that $\mbox{mesh}\left(\rho\right)$ is sufficiently small.
Hence 
\[
\left|\sum_{i=1}^{l}f\left(\nu_{i}\right)\left[g\left(a_{i}\right)-g\left(a_{i-1}\right)\right]-\sum_{i=j}^{m}f\left(\xi_{j}\right)\left[g\left(b_{j}\right)-g\left(b_{j-1}\right)\right]\right|\leq12\cdot2^{-N}S.
\]
provided that $\mbox{mesh}\left(\pi\right)$ and $\mbox{mesh}\left(\rho\right)$
are sufficiently small. Since $N$ may be arbitrary large, we get
the convergence of the approximating sums to an universal limit, which
is the Riemann-Stieltjes integral.\par The estimate (\ref{eq:estimate_integral})
follows directly from the proved convergence of approximating sums
to the Riemann-Stieltjes integral and Lemma \ref{lema}. 
\end{dwd} 
Using Remark \ref{symmetry} and reasoning similarly as in the proof of Theorem \ref{main}, we get the symmetric result.
\begin{theo} \label{main1} Let $f,g:\left[a;b\right]\rightarrow\R$
be two regulated functions which have no common points of discontinuity.
Let $\eta_{0}\geq\eta_{1}\geq\ldots$ and $\theta_{0}\geq\theta_{1}\geq\ldots$
be two sequences of non-negative numbers, such that $\eta_{k}\downarrow0,$
$\theta_{k}\downarrow0$ as $k\rightarrow+\infty.$ Define $\theta_{-1}:=\sup_{a\leq t\leq b}\left|g\left(b\right)-g\left(t\right)\right|$
and 
\[
\tilde{S}:=\sum_{k=0}^{+\infty}2^{k}\theta_{k-1}\cdot\TTV f{\left[a;b\right]}{\eta_{k}}+\sum_{k=0}^{\infty}2^{k}\eta_{k}\cdot\TTV g{\left[a;b\right]}{\theta_{k}}.
\]
If $\tilde{S}<+\infty$ then the Riemann-Stieltjes integral 
$
\int_{a}^{b}f\left(t\right)\mathrm{d}g\left(t\right)
$
exists and one has the following estimate 
\begin{equation}
\left|\int_{a}^{b}f\left(t\right)\mathrm{d}g\left(t\right)-f\left(a\right)\left[g\left(b\right)-g\left(a\right)\right]\right|\leq \tilde{S}.\label{eq:estimate_integral1}
\end{equation}
\end{theo}
From Theorem \ref{main}, Theorem \ref{main1} and Remark \ref{n=1} we also  have
\begin{coro} \label{minST} Let $f,g:\left[a;b\right]\rightarrow\R$
be two regulated functions which have no common points of discontinuity, $\xi \in [a;b]$ and $S$ and $\tilde{S}$ be as in Theorem \ref{main} and Theorem \ref{main1} respectively. If $S<+\ns$ or $\tilde{S}< +\ns$ then the Riemann-Stieltjes integral $\int_{a}^{b}f\left(t\right)\mathrm{d}g\left(t\right)$ exists and one has the following estimate 
\[
\left|\int_{a}^{b}f\left(t\right)\mathrm{d}g\left(t\right)-f\left(\xi\right)\left[g\left(b\right)-g\left(a\right)\right]\right|\leq 
2 \min\{S,\tilde{S} \}.
\]
\end{coro}

\subsection{Young's Theorem and the Lo\'{e}ve-Young inequality}

Recall that for $p>0,$ ${\cal V}^{p}\rbr{[a;b]}$ denotes the family of functions $f:[a;b]\ra \R$ with finite $p-$variation. Note that if $f \in {\cal V}^{p}\rbr{[a;b]}$ then $f$ is regulated. The additional relation we will use, is the following one: if $f\in{\cal V}^{p}\rbr{[a;b]}$ for some
$p\geq 1,$ then for every $\delta>0,$ 
\begin{equation}
\TTV f{\left[a;b\right]}{\delta}\leq V^{p}\left(f,\left[a;b\right]\right)\delta^{1-p}.\label{eq:p_variation}
\end{equation}
As far as we know, the first result of this kind, namely, $\TTV f{\left[a;b\right]}{\delta}\leq C_{f}\delta^{1-p}$
for a continuous function $f:\left[a;b\right]\rightarrow\R$ and some
constant $C_{f}<+\infty$ depending on $f,$ when $f\in{\cal V}^{p}\rbr{[a;b]},$
was proven in \cite[Sect. 6]{TronelVladimirov:2000}. In \cite{TronelVladimirov:2000},
$\TTV f{\left[a;b\right]}{\varepsilon}$ is called $\varepsilon-$variation
and is denoted by $V_{f}(\varepsilon).$ However, being equipped with formula (\ref{TV_def1}), we see that relation (\ref{eq:p_variation})
follows immediately from the inequality: for any $a\leq s<t\leq b,$
\[
\max\left\{ \left|f\left(t\right)-f\left(s\right)\right|-\delta,0\right\} \leq\frac{\left|f\left(t\right)-f\left(s\right)\right|^{p}}{\delta^{p-1}},
\]
which is an obvious consequence of the estimate: $\delta^{p-1}\max\left\{ \left|x\right|-\delta,0\right\} \leq\left|x\right|^{p}$
for any real $x.$ 

Let us denote 
\begin{equation} \label{p_var_norm_def}
\left\Vert f\right\Vert_{p-\text{var},\left[a;b\right]}:=\rbr{V^{p}\left(f,\left[a;b\right]\right)}^{1/p} \text{ and }\left\Vert f\right\Vert _{\text{osc},\left[a;b\right]}:=\sup_{a\leq s<t\leq b}\left|f\left(t\right)-f\left(s\right)\right|.
\end{equation}
Now we are ready to state a Corollary stemming from Theorem \ref{main}, which was one of the main results of \cite{Young:1936}. The second part of this Corollary is a stronger version of the Lo\'{e}ve-Young inequality. By the stronger version we do not mean that we get better constants than in the original inequality (in fact, our constant $C_{p,q}$ is much bigger than the Lo\'{e}ve-Young constant, which is of order $\zeta\rbr{p^{-1}+q^{-1}},$ and it can not be much improved with our methods). We rather mean that the ratio of the Lo\'{e}ve-Young bound to our bound  is always greater than some positive constant, depending on $p$ and $q$ only, but is unbounded from above.
\begin{coro} \label{corol_Young} Let $f,g:\left[a;b\right]\rightarrow\R$
be two functions with no common points of discontinuity.
If $f\in{\cal V}^{p}\rbr{[a;b]}$ and $g\in{\cal V}^{q}\rbr{[a;b]},$ where $p>1,$ $q>1,$
$p^{-1}+q^{-1}>1,$ then the Riemann Stieltjes $\int_{a}^{b}f\left(t\right)\mathrm{d}g\left(t\right)$
exists. Moreover, there exist a constant $C_{p,q},$ depending on $p$ and $q$ only, such that  
\begin{align*}
\left|\int_{a}^{b}f\left(t\right)\mathrm{d}g\left(t\right)-f\left(a\right)\left[g\left(b\right)-g\left(a\right)\right]\right| 
  \leq C_{p,q}\left\Vert f\right\Vert_{p-\emph{var},\left[a;b\right]}^{p-p/q}\left\Vert f\right\Vert _{\emph{osc},\left[a;b\right]}^{1+p/q-p} \left\Vert g\right\Vert_{q-\emph{var},\left[a;b\right]}.
\end{align*}
\end{coro} 
\begin{dwd} By Theorem \ref{main} it is enough to prove that for
some positive sequences $\eta_{-1}=\sup_{a\leq t\leq b}\left|f\left(t\right)-f\left(a\right)\right|,$ $\eta_{0}\geq \eta_1 \geq \ldots$
and $\theta_{0}\geq\theta_{1}\geq\ldots,$ such that $\eta_{k}\downarrow0,$
$\theta_{k}\downarrow0$ as $k\rightarrow+\infty,$ one has 
\begin{align*}
S: & =\sum_{k=0}^{+\infty}2^{k}\eta_{k-1}\cdot\TTV g{\left[a;b\right]}{\theta_{k}}+\sum_{k=0}^{+\infty}2^{k}\theta_{k}\cdot\TTV f{\left[a;b\right]}{\eta_{k}},\\
 & \leq C_{p,q}\left\Vert f\right\Vert_{p-\text{var},\left[a;b\right]}^{p-p/q}\left\Vert f\right\Vert _{\text{osc},\left[a;b\right]}^{1+p/q-p} \left\Vert g\right\Vert_{q-\text{var},\left[a;b\right]}.
\end{align*}
The proof will follow from the proper choice of the sequences $\left(\eta_{k}\right)$
and $\left(\theta_{k}\right).$ Since $p^{-1}+q^{-1}>1,$ we have
$\left(q-1\right)\left(p-1\right)<1.$ Choose \[\alpha=\frac{\sqrt{(q-1)(p-1)}+1}{2}, \quad 
\beta=\sup_{a\leq t\leq b}\left|f\left(t\right)-f\left(a\right)\right|, \quad  \gamma>0
\] and for $k=0,1,\ldots,$ define 
\[
\eta_{k-1}=\beta \cdot  2^{-\left(\alpha^{2}/\left[\left(q-1\right)\left(p-1\right)\right]\right)^{k}+1},
\]
\[
\theta_{k}= \gamma \cdot 2^{-\left(\alpha^{2}/\left[\left(q-1\right)\left(p-1\right)\right]\right)^{k}\alpha/\left(q-1\right)}.
\]
By (\ref{eq:p_variation}) we estimate 
\begin{align*}
\eta_{k-1}\cdot\TTV g{\left[a;b\right]}{\theta_{k}}\leq & \beta \cdot 2^{-\left(\alpha^{2}/\left[\left(q-1\right)\left(p-1\right)\right]\right)^{k}+1}\\
 & \times V^{q}\left(g,\left[a;b\right]\right)\left(\gamma \cdot 2^{-\left(\alpha^{2}/\left[\left(q-1\right)\left(p-1\right)\right]\right)^{k}\alpha/\left(q-1\right)}
\right)^{1-q}\\
= & 2^{-\left(1-\alpha\right)\left(\alpha^{2}/\left[\left(q-1\right)\left(p-1\right)\right]\right)^{k}+1}V^{q}\left(g,\left[a;b\right]\right)\beta\gamma^{1-q},
\end{align*}
and similarly 
\begin{align*}
\theta_{k}\cdot\TTV f{\left[a;b\right]}{\eta_{k}}\leq & \gamma \cdot 2^{-\left(\alpha^{2}/\left[\left(q-1\right)\left(p-1\right)\right]\right)^{k}\alpha/\left(q-1\right)}
\\
 & \times V^{p}\left(f,\left[a;b\right]\right)\left(\beta \cdot 2^{-\left(\alpha^{2}/\left[\left(q-1\right)\left(p-1\right)\right]\right)^{k+1}+1}
\right)^{1-p}\\
= & 2^{-\left(1-\alpha\right)\left(\alpha^{2}/\left[\left(q-1\right)\left(p-1\right)\right]\right)^{k}\alpha/\left(q-1\right)+1-p}V^{p}\left(f,\left[a;b\right]\right)
\beta^{1-p}\gamma.
\end{align*}
Hence 
\begin{eqnarray*}
S & = & \sum_{k=0}^{+\infty}2^{k}\eta_{k-1}\cdot\TTV g{\left[a;b\right]}{\theta_{k}}+\sum_{k=0}^{+\infty}2^{k}\theta_{k}\cdot\TTV f{\left[a;b\right]}{\eta_{k}}\\
 & \leq & \left(\sum_{k=0}^{+\infty}2^{k}2^{-\left(1-\alpha\right)\left(\alpha^{2}/\left[\left(q-1\right)\left(p-1\right)\right]\right)^{k}+1}\right)V^{q}\left(g,\left[a;b\right]\right)\beta\gamma^{1-q}\\
 &  & +\left(\sum_{k=0}^{+\infty}2^{k}2^{-\left(1-\alpha\right)\left(\alpha^{2}/\left[\left(q-1\right)\left(p-1\right)\right]\right)^{k}\alpha/\left(q-1\right)+1-p}\right)V^{p}\left(f,\left[a;b\right]\right)
\beta^{1-p}\gamma.
\end{eqnarray*}
Since $\alpha<1$ and $\alpha^{2}/\left[\left(q-1\right)\left(p-1\right)\right]>1,$
we easily infer that $S<+\infty,$ from which we get that the integral
$\int_{a}^{b}f\left(t\right)\mathrm{d}g\left(t\right)$ exists. 

Moreover,
denoting 
\begin{align*}
C_{p,q} & = \max\left\{ \sum_{k=0}^{+\infty}2^{k+2-\left(1-\alpha\right)\left(\alpha^{2}/\left[\left(q-1\right)\left(p-1\right)\right]\right)^{k}}, \right. \\
& \quad \quad \quad \quad \left.\sum_{k=0}^{+\infty}2^{k+2-\left(1-\alpha\right)\left(\alpha^{2}/\left[\left(q-1\right)\left(p-1\right)\right]\right)^{k}\alpha/\left(q-1\right)-p}\right\} 
\end{align*}
we get 
\begin{align*}
S & \leq\frac{1}{2}C_{p,q}\left(V^{q}\left(g,\left[a;b\right]\right)\beta\gamma^{1-q}+V^{p}\left(f,\left[a;b\right]\right)
\beta^{1-p}\gamma\right).
\end{align*}
Setting in this expression 
$
\gamma=\left(V^{q}\left(g,\left[a;b\right]\right)\right/V^{p}\left(f,\left[a;b\right]\right))^{1/q}
\beta^{p/q}
$
we obtain 
\begin{align*}
S & \leq C_{p,q}\left(V^{q}\left(g,\left[a;b\right]\right)\right)^{1/q}\left(V^{p}\left(f,\left[a;b\right]\right)\right)^{1-1/q}
\beta^{1+p/q-p}
\\ 
& \leq C_{p,q}\left\Vert g\right\Vert_{q-\text{var},\left[a;b\right]}\left\Vert f\right\Vert_{p-\text{var},\left[a;b\right]}^{p-p/q}\left\Vert f\right\Vert _{\text{osc},\left[a;b\right]}^{1+p/q-p}.
\end{align*}
\end{dwd} 
\begin{rem} \label{rem_Young}
Let $f,$ $g,$ $p,$ $q$ and $C_{p,q}$ be the same as in Corollary \ref{corol_Young}. Using Theorem \ref{main1} instead of Theorem \ref{main}, we get the following, similar estimate
\begin{align*}
\left|\int_{a}^{b}f\left(t\right)\mathrm{d}g\left(t\right)-f\left(a\right)\left[g\left(b\right)-g\left(a\right)\right]\right| 
  \leq C_{p,q}\left\Vert f\right\Vert_{p-\emph{var},\left[a;b\right]}\left\Vert g\right\Vert_{q-\emph{var},\left[a;b\right]}^{q-q/p}\left\Vert g\right\Vert _{\emph{osc},\left[a;b\right]}^{1+q/p-q}.
\end{align*}
From Corollary \ref{minST} and the obtained estimates, we also have that for any $\xi \in [a;b]$ 
\begin{align*}
& \left|\int_{a}^{b}f\left(t\right)\mathrm{d}g\left(t\right)-f\left(\xi\right)\left[g\left(b\right)-g\left(a\right)\right]\right| 
 \\
& \leq 2 C_{p,q}
\left\Vert f\right\Vert_{p-\emph{var},\left[a;b\right]} \left\Vert g\right\Vert_{q-\emph{var},\left[a;b\right]}
\min \cbr{\frac{\left\Vert f\right\Vert _{\emph{osc},\left[a;b\right]}^{1+p/q-p}}{\left\Vert f\right\Vert_{p-\emph{var},\left[a;b\right]}^{1+p/q-p} 
}, 
\frac{\left\Vert g\right\Vert _{\emph{osc},\left[a;b\right]}^{1+q/p-q}}{{\left\Vert g\right\Vert_{q-\emph{var},\left[a;b\right]}^{1+q/p-q}}} 
}.
\end{align*}
\end{rem}
\begin{rem} \label{rem_int_V^p}
From Corollary \ref{corol_Young}, reasoning in the similar way as in \cite[p. 456]{Lyons:1994}, we get the following important estimate of the $\Varnormthm {\int_a^{\cdot} f(t) \dd g(t)}{p}$
\[
\Varnormthm {\int_a^{\cdot} f(t) \dd g(t)}{p} \leq \rbr{ C_{p,q}\left\Vert f\right\Vert_{p-\emph{var},\left[a;b\right]}^{p-p/q}\left\Vert f\right\Vert _{\emph{osc},\left[a;b\right]}^{1+p/q-p} + \left\Vert f\right\Vert_{\ns, [a;b]} }\left\Vert g\right\Vert_{q-\emph{var},\left[a;b\right]}.
\]
\end{rem}

\section{The space ${\cal U}^p\rbr{[a;b]}$} \label{space}

Analysing the proof of Corollary \ref{corol_Young} we see that the crucial 
observation we used in the proof were the estimates of the form $\TTV f{[a;b]}{\eta} \leq C \eta^{1-p},$ $\TTV g{[a;b]}{\theta} \leq D \theta^{1-q},$ for any $\eta, \theta >0$ and some constants $C$ and $D.$ These estimates (with $C=V^p\rbr{f,[a;b]}$ and $D=V^q\rbr{g,[a;b]}$) followed from inequality (\ref{eq:p_variation}), but Corollary \ref{corol_Young} remains true (with $\left\Vert f\right\Vert_{p-\text{var},\left[a;b\right]}$ and $\left\Vert g\right\Vert_{q-\text{var},\left[a;b\right]}$ replaced by appropriate constants) as long as $\TTV f{[a;b]}{\eta}$ and $\TTV g{[a;b]}{\theta}$ do not grow faster than $C \eta^{1-p}$ and $D \theta^{1-q}$ when $\eta, \theta \downarrow 0.$ 

For $p\geq 1$ and $f:[a;b] \ra \R$ define 
\begin{equation} \label{eq:TV_p_norm}
\left\Vert f\right\Vert_{p-\text{TV},\left[a;b\right]}:= \sup_{\delta >0} \rbr{\delta^{p-1} \TTV f{[a;b]}\delta}^{1/p}.
\end{equation}
From the very definition of the functional $\left\Vert \cdot \right\Vert_{p-\text{TV},\left[a;b\right]},$ we have that for any $\delta >0,$
\begin{equation} \label{ineq:TV_variation}
\TTV f{[a;b]}\delta \leq \left\Vert f\right\Vert_{p-\text{TV},\left[a;b\right]}^p \delta^{1-p}.
\end{equation}
For $p\geq 1$ let ${\cal U}^p\rbr{[a;b]}$ denote the family of functions $f:[a;b]\ra \R$ such that $\left\Vert f\right\Vert_{p-\text{TV},\left[a;b\right]} < +\ns.$ If $f \in {\cal U}^p\rbr{[a;b]}$ for some $p\geq 1$ then $f$ is regulated. We will prove that the functional defined with formula (\ref{eq:TVnorm}) is a norm and  ${\cal U}^p\rbr{[a;b]}$ equipped with this norm is a Banach space, but first let us state the following counterpart of Corollary \ref{corol_Young}.
\begin{coro} \label{corol_Young1} Let $f,g:\left[a;b\right]\rightarrow\R$
be two functions with no common points of discontinuity.
If $f\in{\cal U}^{p}\rbr{[a;b]}$ and $g\in{\cal U}^{q}\rbr{[a;b]},$ where $p>1,$ $q>1,$
$p^{-1}+q^{-1}>1,$ then the Riemann Stieltjes $\int_{a}^{b}f\left(t\right)\mathrm{d}g\left(t\right)$
exists. Moreover, 
\begin{align*}
\left|\int_{a}^{b}f\left(t\right)\mathrm{d}g\left(t\right)-f\left(a\right)\left[g\left(b\right)-g\left(a\right)\right]\right| 
  \leq C_{p,q} \left\Vert f\right\Vert_{p-\emph{TV},\left[a;b\right]}^{p-p/q}\left\Vert f\right\Vert _{\emph{osc},\left[a;b\right]}^{1+p/q-p} \left\Vert g\right\Vert_{q-\emph{TV},\left[a;b\right]},
\end{align*}
where the constant $C_{p,q}$ is the same as in Corollary \ref{corol_Young}.  
\end{coro} 
The proof of Corollary \ref{corol_Young1} goes along the same lines as the proof of Corollary \ref{corol_Young}, the only difference is that instead of inequality (\ref{eq:p_variation}) one uses (\ref{ineq:TV_variation}). 
We also have the following remark. 
\begin{rem} \label{rem_Young1}
Let $f,$ $g,$ $p,$ $q$ and $C_{p,q}$ be the same as in Corollary \ref{corol_Young1}. Using Theorem \ref{main1} instead of Theorem \ref{main}, we get the following, similar estimate
\begin{align*}
\left|\int_{a}^{b}f\left(t\right)\mathrm{d}g\left(t\right)-f\left(a\right)\left[g\left(b\right)-g\left(a\right)\right]\right| 
  \leq C_{p,q} \left\Vert f\right\Vert_{p-\emph{TV},\left[a;b\right]} \left\Vert g\right\Vert_{q-\emph{TV},\left[a;b\right]}^{q-q/p}\left\Vert g\right\Vert _{\emph{osc},\left[a;b\right]}^{1+q/p-q}.
\end{align*}
From Corollary \ref{minST} and the obtained estimates, we also have that for any $\xi \in [a;b]$ 
\begin{align*}
& \left|\int_{a}^{b}f\left(t\right)\mathrm{d}g\left(t\right)-f\left(\xi\right)\left[g\left(b\right)-g\left(a\right)\right]\right| 
 \\
& \leq 2 C_{p,q}
\left\Vert f\right\Vert_{p-\emph{TV},\left[a;b\right]} 
\left\Vert g\right\Vert_{q-\emph{TV},\left[a;b\right]}
\min \cbr{
\frac{\left\Vert f\right\Vert _{\emph{osc},\left[a;b\right]}^{1+p/q-p}}{\left\Vert f\right\Vert_{p-\emph{TV},\left[a;b\right]}^{1+p/q-p} 
}, 
\frac{\left\Vert g\right\Vert _{\emph{osc},\left[a;b\right]}^{1+q/p-q}}{{\left\Vert g\right\Vert_{q-\emph{TV},\left[a;b\right]}^{1+q/p-q}}}
}.
\end{align*}
\end{rem}
Note that from (\ref{eq:p_variation}) and the definition of $\left\Vert \cdot \right\Vert_{p-\text{TV},\left[a;b\right]}$ (recall formula (\ref{p_var_norm_def})) it follows that $\left\Vert f\right\Vert_{p-\text{TV},\left[a;b\right]} \leq \left\Vert f\right\Vert_{p-\text{var},\left[a;b\right]},$ hence Corollary \ref{corol_Young1} and Remark \ref{rem_Young1} give stronger estimates of the differences \[\left|\int_{a}^{b}f\left(t\right)\mathrm{d}g\left(t\right)-f\left(a\right)\left[g\left(b\right)-g\left(a\right)\right]\right| 
\text{ and } \left|\int_{a}^{b}f\left(t\right)\mathrm{d}g\left(t\right)-f\left(\xi\right)\left[g\left(b\right)-g\left(a\right)\right]\right| 
\]
than Corollary \ref{corol_Young} and Remark \ref{rem_Young}.
\subsection{${\cal U}^{p}\rbr{[a;b]}$ as a Banach space}
Now we are going to prove
\begin{prop} \label{Banach_space}
For any $p\geq1,$ the functional $\left\Vert \cdot \right\Vert_{p-\emph{TV},\left[a;b\right]}$ defined with formula (\ref{eq:TV_p_norm}) is a seminorm and the functional $\left\Vert \cdot \right\Vert_{\text{TV},p,\left[a;b\right]}$ defined with formula (\ref{eq:TVnorm}) is a norm on ${\cal U}^{p}\rbr{[a;b]}.$ ${\cal U}^{p}\rbr{[a;b]}$ equipped with this norm is a Banach space. 
\end{prop}
\begin{dwd}
For $p=1,$ $\left\Vert \cdot \right\Vert_{p-\text{TV},\left[a;b\right]}$ coincides with $V^1\rbr{f,\left[a;b\right]},$ $\left\Vert \cdot \right\Vert_{TV,p,\left[a;b\right]}$ coincides with the $1$-variation norm $\left\Vert f \right\Vert_{var,1,\left[a;b\right]} := \left| f(a) \right|  + V^1\rbr{f,\left[a;b\right]} $ and ${\cal U}^{1}\rbr{[a;b]}$ is simply the same as the space of functions with bounded total variation. Therefore, for the rest of the proof we will assume that $p>1.$

The homogenity of $\left\Vert \cdot \right\Vert_{p-\text{TV},\left[a;b\right]}$  and $\left\Vert \cdot \right\Vert_{TV,p,\left[a;b\right]}$ follows easily from the fact that for $\alpha,\delta>0,$ $\TTV{\alpha f}{\left[a;b\right]}{\alpha\delta}=\alpha\TTV f{\left[a;b\right]}{\delta},$
which is the consequence of the equality \[\left(\left|\alpha f\left(t\right)-\alpha f\left(s\right)\right|-\alpha\delta\right)_{+}=\alpha\left(\left|f\left(t\right)-f\left(s\right)\right|-\delta\right)_{+},\] where $(x)_+$ denotes $\max\cbr{x,0}.$

To prove the triangle inequality, let us take $f,g\in{\cal U}^{p}\rbr{[a;b]}.$ Fix $\varepsilon>0.$
Let $\delta_{0}>0$ and $a\leq t_{0}<t_{1}<\ldots<t_{n}\leq b$ be
such that 
\begin{align}
\left(\delta_{0}^{p-1}\sum_{i=1}^{n}\left(\left|f\left(t_{i}\right)
-f\left(t_{i-1}\right)+g\left(t_{i}\right)-g\left(t_{i-1}\right)\right|-\delta_{0}\right)_{+}\right)^{1/p} & \geq\left\Vert f+g\right\Vert _{p-\text{TV};\left[a;b\right]}-\varepsilon.\label{eq:approx_eps}
\end{align}
By standard calculus, for $x>0$ and $p\geq 1$ we have 
\begin{equation}
\sup_{\delta>0}\delta^{p-1}\left(x-\delta\right)_{+}
=\sup_{\delta\geq0}\delta^{p-1}\left(x-\delta\right)
=c_{p}x^{p}\label{eq:supremum}
\end{equation}
where $c_{p} = (p-1)^{p-1}/p^{p} \in \sbr{2^{-p};1}.$ Denote $x_{0}^{*}=0$
and for $i=1,2,\ldots,n$ define $x_{i}=\left|f\left(t_{i}\right)-f\left(t_{i-1}\right)+g\left(t_{i}\right)-g\left(t_{i-1}\right)\right|.$
Let $x_{1}^{*}\leq x_{2}^{*}\leq\ldots\leq x_{n}^{*}$ be the non-decreasing
re-arrangement of the sequence $\left(x_{i}\right).$ Notice that
by (\ref{eq:supremum}) for $\delta\in\left[x_{j-1}^{*};x_{j}^{*}\right],$
where $j=1,2,\ldots,n,$ one has 
\begin{align*}
& \delta^{p-1}\sum_{i=1}^{n}\left(\left|f\left(t_{i}\right)-f\left(t_{i-1}\right)+g\left(t_{i}\right)-g\left(t_{i-1}\right)\right|-\delta\right)_{+}\\
& =\delta^{p-1}\sum_{i=j}^{n}\left(x_{i}^{*}-\delta\right)=\delta^{p-1}\left(\sum_{i=j}^{n}x_{i}^{*}-\left(n-j+1\right)\delta\right)\\
&=\left(n-j+1\right)\delta^{p-1}\left(\frac{\sum_{i=j}^{n}x_{i}^{*}}{n-j+1}-\delta\right)
\leq\left(n-j+1\right)c_{p}\left(\frac{\sum_{i=j}^{n}x_{i}^{*}}{n-j+1}\right)^{p}.
\end{align*}
Hence 
\begin{align}
& \sup_{\delta>0}\delta^{p-1}\sum_{i=1}^{n}\left(\left|f\left(t_{i}\right)-f\left(t_{i-1}\right)+g\left(t_{i}\right)-g\left(t_{i-1}\right)\right|-\delta\right)_{+}\nonumber \\
& \leq\max_{j=1,2,\dots,n}\left(n-j+1\right)c_{p}\left(\frac{\sum_{i=j}^{n}x_{i}^{*}}{n-j+1}\right)^{p}.\label{eq:nier_1}
\end{align}
On the other hand, 
\begin{align}
&\sup_{\delta>0}\delta^{p-1}\sum_{i=1}^{n}\left(\left|f\left(t_{i}\right)-f\left(t_{i-1}\right)+g\left(t_{i}\right)-g\left(t_{i-1}\right)\right|-\delta\right)_{+}\nonumber \\
&=\sup_{\delta>0}\delta^{p-1}\sum_{i=1}^{n}\left(x_{i}^{*}-\delta\right)_{+}
\geq\sup_{\delta>0}\max_{j=1,2,\dots,n}\delta^{p-1}\sum_{i=j}^{n}\left(x_{i}^{*}-\delta\right)\nonumber \\
& =\max_{j=1,2,\dots,n}\sup_{\delta>0}\delta^{p-1}\sum_{i=j}^{n}\left(x_{i}^{*}-\delta\right) \nonumber \\
& =\max_{j=1,2,\dots,n}\left(n-j+1\right)c_{p}\left(\frac{\sum_{i=j}^{n}x_{i}^{*}}{n-j+1}\right)^{p}.\label{eq:nier_2}
\end{align}
By (\ref{eq:nier_1}) and (\ref{eq:nier_2}) we get 
\begin{align}
&\left(\sup_{\delta>0}\delta^{p-1}\sum_{i=1}^{n}\left(\left|f\left(t_{i}\right)-f\left(t_{i-1}\right)+g\left(t_{i}\right)-g\left(t_{i-1}\right)\right|-\delta\right)_{+}\right)^{1/p} \nonumber \\
&=\max_{j=1,2,\dots,n}\left(n-j+1\right)^{1/p-1}c_{p}^{1/p}\sum_{i=j}^{n}x_{i}^{*}. \label{eq:optim}
\end{align}
Similarly, denoting by $y_{i}^{*}$ and $z_{i}^{*}$ the non-decreasing
rearrangements of the sequences $y_{i}=\left|f\left(t_{i}\right)-f\left(t_{i-1}\right)\right|$
and $z_{i}=\left|g\left(t_{i}\right)-g\left(t_{i-1}\right)\right|$
respectively, we get 
\[
\left(\sup_{\delta>0}\delta^{p-1}\sum_{i=1}^{n}\left(\left|f\left(t_{i}\right)-f\left(t_{i-1}\right)\right|-\delta\right)_{+}\right)^{1/p}=\max_{j=1,2,\dots,n}\left(n-j+1\right)^{1/p-1}c_{p}^{1/p}\sum_{i=j}^{n}y_{i}^{*}
\]
and 
\[
\left(\sup_{\delta>0}\delta^{p-1}\sum_{i=1}^{n}\left(\left|g\left(t_{i}\right)-g\left(t_{i-1}\right)\right|-\delta\right)_{+}\right)^{1/p}=\max_{j=1,2,\dots,n}\left(n-j+1\right)^{1/p-1}c_{p}^{1/p}\sum_{i=j}^{n}z_{i}^{*}
\]
By the triangle inequality and the definition of $y_{i}^{*}$ and
$z_{i}^{*}$ for $j=1,2,\ldots,n,$ we have 
$
\sum_{i=j}^{n}x_{i}^{*}\leq\sum_{i=j}^{n}y_{i}^{*}+\sum_{i=j}^{n}z_{i}^{*}.
$
Hence 
\begin{align*}
&\max_{j=1,2,\dots,n}\left(n-j+1\right)^{1/p-1}c_{p}^{1/p}\sum_{i=j}^{n}x_{i}^{*}\leq \max_{j=1,2,\dots,n}\left(n-j+1\right)^{1/p-1}c_{p}^{1/p}\sum_{i=j}^{n}\left(y_{i}^{*}+z_{i}^{*}\right)\\
& \leq  \max_{j=1,2,\dots,n}\left(n-j+1\right)^{1/p-1}c_{p}^{1/p}\sum_{i=j}^{n}y_{i}^{*} +\max_{j=1,2,\dots,n}\left(n-j+1\right)^{1/p-1}c_{p}^{1/p}\sum_{i=j}^{n}z_{i}^{*}\\
& \leq \left(\sup_{\delta>0}\delta^{p-1}\sum_{i=1}^{n}\left(\left|f\left(t_{i}\right)-f\left(t_{i-1}\right)\right|-\delta\right)_{+}\right)^{1/p}
 \nonumber \\ & \quad +\left(\sup_{\delta>0}\delta^{p-1}\sum_{i=1}^{n}\left(\left|g\left(t_{i}\right)-g\left(t_{i-1}\right)\right|-\delta\right)_{+}\right)^{1/p}
\\
& \leq \left\Vert f\right\Vert _{p-\text{TV},\left[a;b\right]}+\left\Vert g\right\Vert _{p-\text{TV},\left[a;b\right]}.
\end{align*}
Finally, by (\ref{eq:approx_eps}), (\ref{eq:optim}) and the last
estimate, we get 
\[
\left\Vert f+g\right\Vert _{p-\text{TV},\left[a;b\right]}-\varepsilon\leq\left\Vert f\right\Vert _{p-\text{TV},\left[a;b\right]}+\left\Vert g\right\Vert _{p-\text{TV},\left[a;b\right]}.
\]
Sending $\varepsilon$ to $0$ we get the triangle inequality for $\left\Vert \cdot \right\Vert _{p-\text{TV},\left[a;b\right]}.$ The triangle inequality for $\left\Vert \cdot \right\Vert _{TV,p,\left[a;b\right]}$ also holds, since $\left\Vert f \right\Vert _{TV,p,\left[a;b\right]} = \left|f(a)\right|+\left\Vert f \right\Vert _{p-\text{TV},\left[a;b\right]}.$

Now we will prove that the space ${\cal U}^{p}\rbr{[a;b]}$ equipped with the norm $\left\Vert \cdot \right\Vert _{TV,p,\left[a;b\right]}$ is a Banach space. 
To prove this we will need the following inequality 
\begin{equation} \label{TVsplit}
\TTV{f+g}{[a;b]}{\delta_1+\delta_2} \leq  \TTV{f}{[a;b]}{\delta_1}+\TTV{g}{[a;b]}{\delta_2}
\end{equation}
for any $\delta_1,\delta_2 \geq0.$
It follows from the elementary estimate 
\begin{equation} \label{elementary_estimate}
\rbr{|x_1+x_2| - \delta_1 - \delta_2}_+\leq \rbr{|x_1| - \delta_1}_+ + \rbr{|x_2| - \delta_2}_+
\end{equation}
valid for any real $x_1,$ $x_2$ and nonnegative $\delta_1$ and $\delta_2.$
We also have $\TTV f{[a;b]}{\delta} \geq \rbr{\left\Vert f \right\Vert _{\text{osc},\left[a;b\right]} - \delta}_+.$ From this and (\ref{eq:supremum}) it follows that 
\begin{equation} \label{infinite_norm_est}
\left\Vert f \right\Vert _{TV,p,\left[a;b\right]} \geq \left|f(a)\right|+ c_p^{1/p} \left\Vert f \right\Vert _{\text{osc},\left[a;b\right]} \geq \min\cbr{1,c_p^{1/p}} \left\Vert f \right\Vert _{\ns,\left[a;b\right]}.
\end{equation}
Hence any Cauchy sequence $\rbr{f_n}_{n=1}^{\ns}$ in ${\cal U}^{p}\rbr{[a;b]}$ converges uniformly to some $f_{\ns}:[a;b]\ra \R.$ Assume that $\left\Vert f_{\ns} - f_n \right\Vert _{TV,p,\left[a;b\right]} \nrightarrow 0$ as $n\ra +\ns.$ Thus, there exist a positive number $\kappa$, a sequence of positive integers $n_k \ra +\ns$ and a sequence of positive reals $\delta_k,$ $k=1,2,\ldots,$ such that $\delta_k^{p-1} \TTV{f_{n_k} - f_{\ns}}{[a;b]}{\delta_k} \geq \kappa^p.$ Let $N$ be a positive integer such that 
\begin{equation}\label{Cauchy}
\left\Vert f_{m} - f_n \right\Vert _{TV,p,\left[a;b\right]}  < \kappa/2^{1-1/p} \text{ for } m,n\geq N
\end{equation}
and $k_0$ be the minimal positive integer such that $n_{k_0}\geq N.$ For sufficiently large $n\geq N$ we have  $\left\Vert f_{n} - f_{\ns} \right\Vert _{\ns,\left[a;b\right]} \leq \delta_{k_0}/4,$ hence $\Oscnorm {f_{n} - f_{\ns}}\leq \delta_{k_0}/2$ and
\begin{equation} \label{TVzero}
\TTV{f_{n} - f_{\ns}}{[a;b]}{\delta_{k_0}/2} =0.
\end{equation}
Now, by (\ref{TVsplit})
\[
\TTV{f_{n_{k_0}} - f_{\ns}}{[a;b]}{\delta_{k_0}} \leq  \TTV{f_{n_{k_0}} - f_{n}}{[a;b]}{\delta_{k_0}/2}+\TTV{f_{n} - f_{\ns}}{[a;b]}{\delta_{k_0}/2}.
\]
From this and (\ref{TVzero}) we get 
\[
\rbr{\delta_{k_0}/2}^{p-1}\TTV{f_{n_{k_0}} - f_{n}}{[a;b]}{\delta_{k_0}/2} \geq \delta_{k_0}^{p-1} \TTV{f_{n_{k_0}} - f_{\ns}}{[a;b]}{\delta_{k_0}}/2^{p-1} \geq \kappa^p/2^{p-1}
\]
but this (recall (\ref{eq:TVnorm}))  contradicts (\ref{Cauchy}). Thus, the sequence $\rbr{f_n}_{n=1}^{\ns}$ converges in ${\cal U}^{p}\rbr{[a;b]}$ norm  to $f_{\ns}.$ Since the sequence $\rbr{f_n}_{n=1}^{\ns}$ was chosen in an arbitrary way, it proves that ${\cal U}^{p}\rbr{[a;b]}$ is complete. 
\end{dwd}
\begin{rem}
It is easy to see that the space ${\cal U}^{p}\rbr{[a;b]}$ equipped with the norm $\left\Vert \cdot \right\Vert _{TV,p,\left[a;b\right]}$ is not separable. To see this it is enough to consider the family of functions $f_t:[a;b]\ra \cbr{0,1},$ $f_t(s):= {\bf 1}_{\cbr{t}}(s)$ and apply (\ref{infinite_norm_est}). However, we do not know if the subspace of continuous functions in ${\cal U}^{p}\rbr{[a;b]}$ is separable.
\end{rem}
\begin{rem} \label{superadditivity}
From the triangle inequality for $\left\Vert \cdot \right\Vert_{p-\emph{TV},\left[a;b\right]}$ it follows that it is an subadditivie functional of the interval, i.e., for any $p\geq 1,$ $f:[a;b] \ra \R$ and $d \in (a;b),$
\[
\left\Vert f \right\Vert_{p-\emph{TV},\left[a;b\right]} \leq \left\Vert f \right\Vert_{p-\emph{TV},\left[a;d\right]}  + \left\Vert f \right\Vert_{p-\emph{TV},\left[d;b\right]}. 
\]
To see this it is enough to consider the following decomposition $f(t) = f_1(t)+f_2(t),$ $f_1(t) = f(t) {\bf 1}_{[a;d]}(t) + f(d){\bf 1}_{(d;b]}(t),$ $f_2(t) = f(t) {\bf 1}_{(d;b]}(t) - f(d){\bf 1}_{(d;b]}(t).$ We naturally have 
\begin{align*}
&\left\Vert f \right\Vert_{p-\emph{TV},\left[a;b\right]} = \left\Vert f_1 + f_2 \right\Vert_{p-\emph{TV},\left[a;b\right]} \\
&\leq \left\Vert f_1 \right\Vert_{p-\emph{TV},\left[a;b\right]} + \left\Vert f_2 \right\Vert_{p-\emph{TV},\left[a;b\right]}  \\
& = \left\Vert f \right\Vert_{p-\emph{TV},\left[a;d\right]} + \left\Vert f \right\Vert_{p-\emph{TV},\left[d;b\right]}.
\end{align*}

However, for $p>1$ we have no longer the superadditivity of $\left\Vert \cdot \right\Vert_{p-\emph{TV},\left[a;b\right]}^p,$ which holds for $\left\Vert \cdot \right\Vert_{p-\emph{var},\left[a;b\right]}^p = V^p\rbr{\cdot,[a;b]}.$ To see this it is enough to consider for $x> p/(p-1)$ the function $f_x:[-1;1]\ra \cbr{-1,0,1},$ $f_x(t) = {\bf 1}_{(-1;1)}(t)-(x-1){\bf 1}_{\cbr{1}}(t).$ Since the function $\delta\mapsto\delta^{p-1}\left(x-\delta\right)_{+}$
attains its global maximum at the unique point $\delta_x=x\left(p-1\right)/p,$
for any $x> p/\left(p-1\right)$ we have $\delta_x>1$ and
\begin{align*}
&\left\Vert f_x\right\Vert _{p-\emph{TV};\left[-1;1\right]}^{p}=\sup_{\delta>0}\delta^{p-1}\left(\left(1-\delta\right)_{+}+\left(x-\delta\right)_{+}\right)\\
&=\max\left\{ \sup_{0<\delta<1}\delta^{p-1}\left(\left(1-\delta\right)_{+}+\left(x-\delta\right)_{+}\right),\sup_{\delta\geq1}\delta^{p-1}\left(\left(1-\delta\right)_{+}+\left(x-\delta\right)_{+}\right)\right\} \\
&=\max\left\{ \sup_{0<\delta<1}\delta^{p-1}\left(\left(1-\delta\right)_{+}+\left(x-\delta\right)_{+}\right),\sup_{\delta\geq1}\delta^{p-1}\left(x-\delta\right)_{+}\right\} \\
&<\sup_{\delta>0}\delta^{p-1}\left(1-\delta\right)_{+}+\delta_{x}^{p-1}\left(x-\delta_{x}\right)_{+}=\left\Vert f_x\right\Vert _{p-\emph{TV};\left[-1;0\right]}^{p}+\left\Vert f_x\right\Vert _{p-\emph{TV};\left[0;1\right]}^{p}.
\end{align*}
\end{rem}
\begin{rem} \label{rem:no_splitting}
Opposite than for the seminorm $\left\Vert \cdot\right\Vert _{p-\emph{var};\left[a;b\right]}$
(cf. \cite[Proposition 5.8]{Friz:2010fk}) for each $p>1$ there exist continuous
functions $f\in{\cal U}^{p}\left(\left[a;b\right]\right)$ such that
for any $c\in\left(a;b\right),$ $\left\Vert f\right\Vert _{p-\emph{TV};\left[a;c\right]}\geq1.$
To see this, for $p>1$ it is enough to consider the following continuous function $\varphi\in{\cal U}^{p}\left(\left[0;1\right]\right):$ on the interval
$\left[2^{-n};2^{-n+1}\right],$ $n=1,2,\ldots,$ the function $\varphi$
has $\left\lceil 2^{np-1}\right\rceil $ ``zigzags'' of magnitude
$2^{-n+1}.$ More precisely, on the interval $\left[2^{-n};2^{-n+1}\right],$
$n=1,2,\ldots,$ it attains $\alpha\left(n\right)=\left\lceil 2^{np-1}\right\rceil $
times the value $2^{-n+1}$ at some points $2^{-n}<t_{1}^{\left(n\right)}<t_{2}^{\left(n\right)}<\ldots<t_{\alpha\left(n\right)}^{\left(n\right)}<2^{-n+1}$
and value $0$ at some points $s_{i}^{\left(n\right)}\in\left(t_{i}^{\left(n\right)};t_{i+1}^{\left(n\right)}\right),$
$i=1,2,\ldots,\alpha\left(n\right)-1$ and at the points $s_{0}^{\left(n\right)}=2^{-n}$
and $s_{\alpha\left(n\right)}^{\left(n\right)}=2^{-n+1}.$ Moreover,
$\varphi$ is linear on each interval of the form $\left[s_{i-1}^{\left(n\right)};t_{i}^{\left(n\right)}\right]$
and $\left[t_{i}^{\left(n\right)};s_{i}^{\left(n\right)}\right],$
$i=1,2,\ldots,\alpha\left(n\right)-1.$ Denoting $\delta\left(n\right)=2^{-n}$
we easily calculate 
\[
\left(\delta\left(n\right)\right)^{p-1}TV^{\delta\left(n\right)}\left(\varphi,\left[2^{-n};2^{-n+1}\right]\right)\geq\left(2^{-n}\right)^{p-1}2\cdot 2^{np-1}2^{-n}=1,
\]
hence $\left\Vert \varphi\right\Vert _{p-\emph{TV};\left[2^{-n};2^{-n+1}\right]}\geq1.$
On the other hand, $\varphi\in{\cal U}^{p}\left(\left[0;1\right]\right)$
since for $\delta\in\left[2^{-n};2^{-n+1}\right),$ $n=1,2,\ldots,$
we have 
\begin{align*}
TV^{\delta}\left(\varphi,\left[0;1\right]\right) & \leq TV\left(\varphi,\left[2^{-n};1\right]\right)\leq\sum_{k=1}^{n}TV\left(\varphi,\left[2^{-k};2^{-k+1}\right]\right)\\
 & =\sum_{k=1}^{n}2\left(2^{kp-1}+1\right)2^{-k+1}\leq\frac{4}{2^{p-1}-1}2^{\left(n+1\right)\left(p-1\right)}
\end{align*}
hence 
\[
\sup_{\delta>0}\delta^{p-1}TV^{\delta}\left(\varphi,\left[0;1\right]\right)\leq\frac{4}{2^{p-1}-1}\sup_{n=1,2,\ldots}2^{\left(-n+1\right)\left(p-1\right)}2^{\left(n+1\right)\left(p-1\right)}=\frac{4\cdot2^{2\left(p-1\right)}}{2^{p-1}-1}.
\]
\end{rem}
\subsection{Relationships between the spaces ${\cal U}^{p}\rbr{[a;b]}$ and ${\cal V}^{p}\rbr{[a;b]}$}
From inequality (\ref{eq:p_variation}) it follows that for $p\geq 1$
\[
\left\Vert f \right\Vert_{TV,p,\left[a;b\right]} \leq \left\Vert f \right\Vert_{var,p,\left[a;b\right]},
\]
where $\left\Vert f \right\Vert_{var,p,\left[a;b\right]} := \left| f(a) \right|  + \left\Vert f \right\Vert_{p-\text{var},\left[a;b\right]} = \left| f(a) \right| + \rbr{V^p\rbr{f,[a;b]}}^{1/p},$ hence
${\cal V}^{p}\rbr{[a;b]} \subset {\cal U}^{p}\rbr{[a;b]}$ for $p\geq 1.$ This result (proven in a different way) is due to Tronel and Vladimirov, see \cite[Theorem 17]{TronelVladimirov:2000}. They also show, constructing a simple example of a 'jigsaw' function $f$ such that $f\in {\cal U}^{p}\rbr{[a;b]}$ but $f\notin {\cal V}^{p}\rbr{[a;b]}$, see \cite[pp. 95--96]{TronelVladimirov:2000}, that for $p>1$ we have the strict inclusion.  
Another interesting examples, for $p=2,$ come from the theory of stochastic processes. For any $T>0,$ almost every path of a standard Brownian motion $B$ or a continuous semimartingale $X$ with positive quadratic variation, $\langle X\rangle_T > 0,$ is not an element of ${\cal V}^2\rbr{[0;T]}$ (see  \cite{Levy:1940} for the Brownian motion case, the more general, semimartingale case follows from this and the Dumbis, Dubins-Schwarz Theorem, see e.g. \cite[Chap. V, Theorem 1.6]{RevuzYor:2005}). However, from \cite[Theorem 1]{LochowskiMilosSPA:2013} we get that $B \in {\cal U}^2\rbr{[0;T]}$ and $X \in {\cal U}^2\rbr{[0;T]}$ almost surely. 

Tronel and Vladimirov prove also that for any $q>p\geq 1,$ ${\cal U}^{p}\rbr{[a;b]} \subset {\cal V}^{q}\rbr{[a;b]}.$ Using formula (\ref{TV_def1}) we obtain the following quantitative result. 
\begin{prop} \label{prop:tv_p_norm_var_q_norm}
For any $q>p\geq 1$ and $f\in {\cal U}^{p}\rbr{[a;b]}$ we have $f\in {\cal V}^{q}\rbr{[a;b]}$ and 
\begin{equation}\label{tv_var_rel_p_q}
\left\Vert f \right\Vert_{q-\emph{var},\left[a;b\right]} \leq \rbr{\frac{2^{q+p-1}}{2^{q-p}-1}}^{1/q} \left\Vert f \right\Vert_{\emph{osc},\left[a;b\right]}^{1-p/q} \left\Vert f \right\Vert_{p-\emph{TV},\left[a;b\right]}^{p/q}.
\end{equation}
\end{prop}
\begin{dwd}
The first part of the assertion follows from inequality (\ref{tv_var_rel_p_q}). To prove this inequality consider
the partition $\pi=\left\{ a\leq t_{0}<t_{1}<\ldots<t_{n}\leq b\right\} .$
For $j=1,2,\ldots$ let us define 
\[
I_{j}=\left\{ i\in\left\{ 1,2,\ldots,n\right\} :\left|f\left(t_{i}\right)-f\left(t_{i-1}\right)\right|\in\left[2^{-j}\left\Vert f\right\Vert _{\text{osc};\left[a;b\right]};2^{-j+1}\left\Vert f\right\Vert _{\text{osc};\left[a;b\right]}\right]\right\} 
\]
and $\delta\left(j\right):=2^{-j-1}\left\Vert f\right\Vert _{\text{osc};\left[a;b\right]}.$
Naturally, for $i\in I_{j},$ 
\[
\left|f\left(t_{i}\right)-f\left(t_{i-1}\right)\right|-\delta\left(j\right)\geq\frac{1}{2}\left|f\left(t_{i}\right)-f\left(t_{i-1}\right)\right|
\]
and since $\left\{ 1,2,\ldots,n\right\} =\bigcup_{j=1}^{+\infty}I_{j},$
for $q>p$ we estimate 
\begin{align*}
\sum_{i=1}^{n}\left|f\left(t_{i}\right)-f\left(t_{i-1}\right)\right|^{q} & =\sum_{j=1}^{+\infty}\sum_{i\in I_{j}}\left|f\left(t_{i}\right)-f\left(t_{i-1}\right)\right|^{q}\\
 & \leq\sum_{j=1}^{+\infty}\left(2^{-j+1}\left\Vert f\right\Vert _{\text{osc};\left[a;b\right]}\right)^{q-1}\sum_{i\in I_{j}}\left|f\left(t_{i}\right)-f\left(t_{i-1}\right)\right|\\
 & \leq\sum_{j=1}^{+\infty}\left(2^{-j+1}\left\Vert f\right\Vert _{\text{osc};\left[a;b\right]}\right)^{q-1}2\sum_{i\in I_{j}}\max\left\{ \left|f\left(t_{i}\right)-f\left(t_{i-1}\right)\right|-\delta\left(j\right),0\right\} \\
 & \leq\sum_{j=1}^{+\infty}\left(2^{-j+1}\left\Vert f\right\Vert _{\text{osc};\left[a;b\right]}\right)^{q-1}2\TTV f{\left[a;b\right]}{\delta\left(j\right)}\\
 & \leq\sum_{j=1}^{+\infty}\left(2^{-j+1}\left\Vert f\right\Vert _{\text{osc};\left[a;b\right]}\right)^{q-1}2\left\Vert f \right\Vert_{p-\text{TV},\left[a;b\right]}^{p}\left(2^{-j-1}\left\Vert f\right\Vert _{\text{osc};\left[a;b\right]}\right)^{1-p}\\
 & =\frac{2^{q+p-1}}{2^{q-p}-1}\left\Vert f\right\Vert _{\text{osc};\left[a;b\right]}^{q-p}\left\Vert f \right\Vert_{p-\text{TV},\left[a;b\right]}^{p}.
\end{align*}
Since the partition $\pi$ was chosen in an arbitrary way, we get (\ref{tv_var_rel_p_q}).
\end{dwd}
\begin{rem}
The relationship between $\TVnormthm {\cdot}{q}$ and $\TVnormthm {\cdot}{p}$ norms is simpler. For $\delta > \Oscnormthm{f},$ $\rbr{\left| f(t) - f(s) \right| - \delta}_+=0$ hence for $q > p \geq 1$ and any $\delta >0$ \[\delta^{q-1} \rbr{\left| f(t) - f(s) \right| - \delta}_+ \leq {\Oscnormthm{f}^{q-p}} \delta^{p-1} \rbr{\left| f(t) - f(s) \right| - \delta}_+.\] Thus we have  
\[
\TVnormthm {f}{q} \leq {\Oscnormthm{f}^{1-p/q}} \TVnormthm {f}{q}^{p/q}. 
\]
From the estimate $\TTV f{[a;b]}{\delta} \geq \rbr{\left\Vert f \right\Vert _{\text{osc},\left[a;b\right]} - \delta}_+$ and (\ref{eq:supremum}) it follows also that for $p\geq 1$
\begin{equation} \label{ineq_TV_osc}
\TVnormthm{f}{p} \geq c_p^{1/p} \Oscnormthm{f} \geq 2^{-1} \Oscnormthm{f}.
\end{equation}
\end{rem}

\subsection{Integration of functions from ${\cal U}^{p}\rbr{[a;b]}$}
In \cite[Section 2]{Lyons:1994} there are considered $\left\Vert \cdot \right\Vert_{p-\text{var},\left[a;b\right]}$ norms of the integrals of the form $[a;b] \ni t \mapsto \int_a^t f(s) \dd g(s),$ with $f \in {\cal V}^{p}\rbr{[a;b]}$ and $g \in {\cal V}^{q}\rbr{[a;b]},$ where $p>1,$ $q>1$ and $p^{-1}+q^{-1}>1.$   Now, we turn to investigate the $\left\Vert \cdot \right\Vert_{p-\text{TV},\left[a;b\right]}$ norms of similar integrals, but with $f \in {\cal U}^{p}\rbr{[a;b]}$ and $g \in {\cal U}^{q}\rbr{[a;b]}.$ 
\begin{theo} \label{thm_integral}
Assume that $f \in {\cal U}^{p}\rbr{[a;b]}$ and $g \in {\cal U}^{q}\rbr{[a;b]}$ for some $p>1,$ $q>1,$ such that $p^{-1}+q^{-1}>1$ and they have no common points of discontinuity. Then there
exist a constant $D_{p,q}<+\infty,$ depending on $p$ and $q$ only,
such that 
\[ \left\Vert \int_{a}^{\cdot}\left[f\left(s\right)-f\left(a\right)\right]\mathrm{d}g\left(s\right)  \right\Vert_{q-\emph{TV},\left[a;b\right]} 
\leq D_{p,q}\left\Vert f\right\Vert _{p-\emph{TV},\left[a;b\right]}^{p-p/q}\left\Vert f\right\Vert _{\emph{osc},\left[a;b\right]}^{1+p/q-p}\left\Vert g\right\Vert _{q-\emph{TV},\left[a;b\right]}.
\]
\end{theo}
In this case we have no longer the supperadditivity property of the functional $\left\Vert \cdot \right\Vert_{p-\text{TV},\left[a;b\right]}^p$ as the function of  interval (see Remark \ref{superadditivity}), hence the method of the proof of Theorem \ref{thm_integral} will be different than the proofs of related estimates in \cite{Lyons:1994}. It will be similar to the proof of Corollary \ref{corol_Young}. 
We will need the following lemma.
\begin{lema} \label{lema1}
Let $f,g:\left[a;b\right]\rightarrow\R$
be two regulated functions which have no common points of discontinuity 
and $\delta_{0}\geq\delta_{1}\geq\ldots,$ $\varepsilon_{0}\geq\varepsilon_{1}\geq\ldots$
be two sequences of non-negative numbers, such that $\delta_{k}\downarrow0,$
$\varepsilon_{k}\downarrow0$ as $k\rightarrow+\infty.$ Assume that for $\delta_{-1}:=\sup_{a\leq t\leq b}\left|f\left(t\right)-f\left(a\right)\right|$
and 
\[
S=\sum_{k=0}^{+\infty}2^{k}\delta_{k-1}\cdot\TTV g{\left[a;b\right]}{\varepsilon_{k}}+\sum_{k=0}^{\infty}2^{k}\varepsilon_{k}\cdot\TTV f{\left[a;b\right]}{\delta_{k}}
\]
we have $S < +\ns.$
Defining 
\[
\gamma:=2\sum_{k=0}^{+\infty}2^{k}\varepsilon_{k}\cdot\TTV f{\left[a;b\right]}{\delta_{k}}
\]
we get 
\begin{align*}
\thmTTV{\int_{a}^{\cdot}\left[f\left(s\right)-f\left(a\right)\right]\mathrm{d}g
\left(s\right)}{\left[a;b\right]}{\gamma} & \leq\sum_{k=0}^{+\infty}2^{k}\delta_{k-1}\cdot\TTV g{\left[a;b\right]}{\varepsilon_{k}}.
\end{align*}
\end{lema}
\begin{dwd} We proceed similarly as in the proof of Lemma \ref{lema}.
Define $g_{0}=g,$ $f_{0}=f,$ $g_{1}:=g_{0}-g_{0}^{\varepsilon_{0}},$
$f_{1}:=f_{0}-f_{0}^{\delta_{0}},$ where $g_{0}^{\varepsilon_{0}}$
is regulated and such that 
\[
\left\Vert g_{0}-g_{0}^{\varepsilon_{0}}\right\Vert _{\infty}\leq\frac{1}{2}\varepsilon_{0}\mbox{ and }\TTV{g_{0}^{\varepsilon_{0}}}{\left[a;b\right]}0=\TTV{g_{0}}{\left[a;b\right]}{\varepsilon_{0}},
\]
$f_{0}^{\delta_{0}}$ is regulated and such that 
\[
\left\Vert f_{0}-f_{0}^{\delta_{0}}\right\Vert _{\infty}\leq\frac{1}{2}\delta_{0}\mbox{ and }\TTV{f_{0}^{\delta_{0}}}{\left[a;b\right]}0=\TTV{f_{0}}{\left[a;b\right]}{\delta_{0}},
\]
and for $k=2,3,\ldots,$ $g_{k}:=g_{k-1}-g_{k-1}^{\varepsilon_{k-1}},$
$f_{k}:=f_{k-1}-f_{k-1}^{\delta_{k-1}}$ are defined similarly as
$g_{1}$ and $f_{1}.$ By the linearity of the RSI with respect
to the integrator, integrating by parts, for $t\in\left[a;b\right],$
$r=1,2,\ldots,$ we have 
\begin{align}
\int_{a}^{t}\left[f\left(s\right)-f\left(a\right)\right]\mathrm{d}g\left(s\right) & =\int_{a}^{t}\left[f_{0}\left(s\right)-f_{0}\left(a\right)\right]\mathrm{d}g^{\varepsilon_{0}}\left(s\right)+\int_{a}^{t}\left[f_{0}\left(s\right)-f_{0}\left(a\right)\right]\mathrm{d}g_{1}\left(s\right)\nonumber \\
 & =\int_{a}^{t}\left[f_{0}\left(s\right)-f_{0}\left(a\right)\right]\mathrm{d}g_{0}^{\varepsilon_{0}}\left(s\right)+\int_{a}^{t}\left[g_{1}\left(t\right)-g_{1}\left(s\right)\right]\mathrm{d}f_{0}\left(s\right)\nonumber \\
 & =\int_{a}^{t}\left[f_{0}\left(s\right)-f_{0}\left(a\right)\right]\mathrm{d}g_{0}^{\varepsilon_{0}}\left(s\right)+\int_{a}^{t}\left[g_{1}\left(t\right)-g_{1}\left(s\right)\right]\mathrm{d}f_{0}^{\delta_{0}}\left(s\right)\nonumber \\
 & \quad +\int_{a}^{t}\left[g_{1}\left(t\right)-g_{1}\left(s\right)\right]\mathrm{d}f_{1}\left(s\right)\nonumber \\
 & =\int_{a}^{t}\left[f_{0}\left(s\right)-f_{0}\left(a\right)\right]\mathrm{d}g_{0}^{\varepsilon_{0}}\left(s\right)+\int_{a}^{t}\left[g_{1}\left(t\right)-g_{1}\left(s\right)\right]\mathrm{d}f_{0}^{\delta_{0}}\left(s\right)\nonumber \\
 & \quad +\int_{a}^{t}\left[f_{1}\left(s\right)-f_{1}\left(a\right)\right]\mathrm{d}g_{1}\left(s\right)=\ldots\nonumber \\
 & =\sum_{k=0}^{r-1}\left(\int_{a}^{t}\left[f_{k}\left(s\right)-f_{k}\left(a\right)\right]\mathrm{d}g_{k}^{\varepsilon_{k}}\left(s\right)+\int_{a}^{t}\left[g_{k+1}\left(t\right)-g_{k+1}\left(s\right)\right]\mathrm{d}f_{k}^{\delta_{k}}\left(s\right)\right)\nonumber \\
 & \quad +\int_{a}^{t}\left[f_{r}\left(s\right)-f_{r}\left(a\right)\right]\mathrm{d}g_{r}\left(s\right).\label{eq:summaa}
\end{align}
By Theorem \ref{main}, we easily estimate that 
\begin{align}
\left|\int_{a}^{t}\left[f_{r}\left(s\right)-f_{r}\left(a\right)\right]\mathrm{d}g_{j}\left(s\right)\right|\leq & \sum_{k=r}^{+\infty}2^{k}\delta_{k-1}\cdot\TTV g{\left[a;t\right]}{\varepsilon_{k}}+\sum_{k=r}^{\infty}2^{k}\varepsilon_{k}\cdot\TTV f{\left[a;t\right]}{\delta_{k}}\label{eq:tv_gamma_est1}
\end{align}
for $r=1,2,\ldots.$ 
Moreover, for $k=0,1,\ldots,$ similarly as in the proof of Lemma
\ref{lema}, we estimate 
\begin{equation}
\left|\int_{a}^{t}\left[g_{k+1}\left(t\right)-g_{k+1}\left(s\right)\right]\mathrm{d}f_{k}^{\delta_{k}}\left(s\right)\right|\leq\varepsilon_{k}^{k}\TTV{f_{k}^{\delta_{k}}}{\left[a;t\right]}0\leq2\varepsilon_{k}^{k}\TTV f{\left[a;t\right]}{\delta_{k}},\label{eq:tv_gamma_est2}
\end{equation}
and 
\begin{align}
\TTV{\int_{a}^{\cdot}\left[f_{k}\left(s\right)-f_{k}\left(a\right)\right]\mathrm{d}g_{k}^{\varepsilon_{k}}\left(s\right)}{\left[a;b\right]}0 & \leq\delta_{k-1}\TTV{g_{k}^{\varepsilon_{k}}}{\left[a;b\right]}0\nonumber \\
 & \leq2^{k}\delta_{k-1}\cdot\TTV g{\left[a;b\right]}{\varepsilon_{k}}.\label{eq:tv_gamma_est3}
\end{align}
(Notice that for the function $F_k(t) :=  \int_{a}^{t}\left[g_{k+1}\left(t\right)-g_{k+1}\left(s\right)\right]\mathrm{d}f_{k}^{\delta_{k}}\left(s\right)$
we could not obtain a similar estimate as (\ref{eq:tv_gamma_est3}). This is due to the fact that $F_k(t_2) - F_k(t_1)$ can not be expressed as the integral $\int_{t_1}^{t_2}\left[g_{k+1}\left(t\right)-g_{k+1}\left(s\right)\right]\mathrm{d}f_{k}^{\delta_{k}}\left(s\right)$
.)
Defining 
\[
\gamma\left(r\right):=2\sum_{k=0}^{r-1}2^{k}\varepsilon_{k}\cdot\TTV f{\left[a;b\right]}{\delta_{k}}+2\sum_{k=r}^{+\infty}2^{k}\delta_{k-1}\cdot\TTV g{\left[a;b\right]}{\varepsilon_{k}},
\]
from (\ref{eq:summaa}), (\ref{eq:tv_gamma_est2}) and (\ref{eq:tv_gamma_est1})
we get 
\begin{align*}
 & \left|\int_{a}^{t}\left[f\left(s\right)-f\left(a\right)\right]\mathrm{d}g\left(s\right)-\sum_{k=0}^{r-1}\int_{a}^{t}\left[f_{k}\left(s\right)-f_{k}\left(a\right)\right]\mathrm{d}g_{k}^{\varepsilon_{k}}\left(s\right)\right|\\
 & \leq\sum_{k=0}^{r-1}\left|\int_{a}^{t}\left[g_{k+1}\left(t\right)-g_{k+1}\left(s\right)\right]\mathrm{d}f_{k}^{\delta_{k}}\left(s\right)\right|+\left|\int_{a}^{t}\left[f_{r}\left(s\right)-f_{r}\left(a\right)\right]\mathrm{d}g_{j}\left(s\right)\right|\\
 & \leq\frac{1}{2}\gamma\left(r\right)
\end{align*}
for any $t\in\left[a;b\right].$ Let us notice that by the very definition
of the truncated variation, $$\TTV{\int_{a}^{\cdot}\left[f\left(s\right)-f\left(a\right)\right]\mathrm{d}g
\left(s\right)}{\left[a;t\right]}{\gamma} $$ is bouned from above
by the variation of any function approximating $\int_{a}^{\cdot}f\left(s\right)\mathrm{d}g\left(s\right)$
with accuracy $\gamma/2.$ By this variational property of the truncated
variation and by (\ref{eq:tv_gamma_est3}) we get 
\begin{align*}
\TTV{\int_{a}^{\cdot}\left[f\left(s\right)-f\left(a\right)\right]\mathrm{d}g\left(s\right)}{\left[a;b\right]}{\gamma\left(r\right)}\leq & \TTV{\sum_{k=0}^{r-1}\int_{a}^{\cdot}\left[f_{k}\left(s\right)-f_{k}\left(a\right)\right]\mathrm{d}g_{k}^{\varepsilon_{k}}\left(s\right)}{\left[a;b\right]}0\\
\leq & \sum_{k=0}^{r-1}2^{k}\delta_{k-1}\cdot\TTV g{\left[a;b\right]}{\varepsilon_{k}}.
\end{align*}
Proceeding with $r$ to $+\infty$ we get the assertion.
\end{dwd}
Now we are ready to prove Theorem \ref{thm_integral}.
\begin{dwd} 
Let
$\gamma>0.$ We choose $\alpha=\frac{\sqrt{(q-1)(p-1)}+1}{2}<1,$  set 
\[
\beta:=\frac{\gamma}{\left(\sum_{k=0}^{+\infty}2^{k+2-\left(1-\alpha\right)\left(\alpha^{2}/\left[\left(q-1\right)\left(p-1\right)\right]\right)^{k}\alpha/\left(q-1\right)-p}\right)\left\Vert f \right\Vert_{p-\text{TV},\left[a;b\right]}^p\left\Vert f\right\Vert _{\text{osc},\left[a;b\right]}^{1/\left(q-1\right)+1-p}}
\]
and for $k=0,1,\ldots,$ define 
\[
\delta_{k-1}=2^{-\left(\alpha^{2}/\left[\left(q-1\right)\left(p-1\right)\right]\right)^{k}+1}\sup_{a\leq t\leq b}\left|f\left(t\right)-f\left(a\right)\right|,
\]
\[
\varepsilon_{k}=2^{-\left(\alpha^{2}/\left[\left(q-1\right)\left(p-1\right)\right]\right)^{k}\alpha/\left(q-1\right)}\sup_{a\leq t\leq b}\left|f\left(t\right)-f\left(a\right)\right|^{1/\left(q-1\right)}\beta.
\]
Using (\ref{ineq:TV_variation}), similarly as in the proof of Corollary
\ref{corol_Young} { we estimate 
\begin{align*}
\sum_{k=0}^{+\infty}2^{k}\delta_{k-1}\cdot\TTV g{\left[a;b\right]}{\varepsilon_{k}} & \leq\left(\sum_{k=0}^{+\infty}2^{k+1-\left(1-\alpha\right)\left(\alpha^{2}/\left[\left(q-1\right)\left(p-1\right)\right]\right)^{k}}\right)\left\Vert g \right\Vert_{q-\text{TV},\left[a;b\right]}^q\beta^{1-q}
\end{align*}
and 
\begin{align*}
\tilde{\gamma}: & =2\sum_{k=0}^{+\infty}2^{k}\varepsilon_{k}\cdot\TTV f{\left[a;b\right]}{\delta_{k}}\\
 & \leq\left(\sum_{k=0}^{+\infty}2^{k+2-\left(1-\alpha\right)\left(\alpha^{2}/\left[\left(q-1\right)\left(p-1\right)\right]\right)^{k}\alpha/\left(q-1\right)-p}\right)
\left\Vert f \right\Vert_{p-\text{TV},\left[a;b\right]}^p
\left\Vert f\right\Vert_{\text{osc},\left[a;b\right]}^{1/\left(q-1\right)+1-p}\beta\\
 & =\gamma.
\end{align*}
By the monotonicity of the truncated variation, Lemma \ref{lema1} and the
last two estimates we get 
\begin{align*}
 & \TTV{\int_{a}^{\cdot}\left[f\left(s\right)-f\left(a\right)\right]\mathrm{d}g\left(s\right)}{\left[a;b\right]}{\gamma}\leq\TTV{\int_{a}^{\cdot}\left[f\left(s\right)-f\left(a\right)\right]\mathrm{d}g\left(s\right)}{\left[a;b\right]}{\tilde{\gamma}}\\
 & \leq\sum_{k=0}^{+\infty}2^{k}\delta_{k-1}\cdot\TTV g{\left[a;b\right]}{\varepsilon_{k}}\\
 & \leq\left(\sum_{k=0}^{+\infty}2^{k+1-\left(1-\alpha\right)\left(\alpha^{2}/\left[\left(q-1\right)\left(p-1\right)\right]\right)^{k}}\right)
\left\Vert g \right\Vert_{q-\text{TV},\left[a;b\right]}^q\beta^{1-q}\\
 & =\tilde{D}_{p,q}\left\Vert f \right\Vert_{p-\text{TV},\left[a;b\right]}^{pq-p} \left\Vert f\right\Vert _{\text{osc},\left[a;b\right]}^{p+q-pq}\left\Vert g \right\Vert_{q-\text{TV},\left[a;b\right]}^q \gamma^{1-q},
\end{align*}
where 
\begin{align*}
\tilde{D}_{p,q}= & \left(\sum_{k=0}^{+\infty}2^{k+1-\left(1-\alpha\right)\left(\alpha^{2}/\left[\left(q-1\right)\left(p-1\right)\right]\right)^{k}}\right)\\
 & \times\left(\sum_{k=0}^{+\infty}2^{k+2-\left(1-\alpha\right)\left(\alpha^{2}/\left[\left(q-1\right)\left(p-1\right)\right]\right)^{k}\alpha/\left(q-1\right)-p}\right)^{q-1}.
\end{align*}
}
From this and the definition of $\left\Vert \cdot \right\Vert_{q-\text{TV},\left[a;b\right]} $ we have 
\[
\left\Vert \int_{a}^{\cdot}\left[f\left(s\right)-f\left(a\right)\right]\mathrm{d}g\left(s\right) \right\Vert_{q-\text{TV},\left[a;b\right]} \leq D_{p,q}\left\Vert f \right\Vert_{p-\text{TV},\left[a;b\right]}^{p-p/q} \left\Vert f\right\Vert _{\text{osc},\left[a;b\right]}^{p/q+1-p}\left\Vert g \right\Vert_{q-\text{TV},\left[a;b\right]},
\]
where $D_{p,q} = \tilde{D}_{p,q}^{1/q}.$
\end{dwd}
Theorem \ref{thm_integral} and estimate (\ref{ineq_TV_osc}) imply 
\begin{coro} \label{coro_Love_Young_TV}
Assume that $f \in {\cal U}^{p}\rbr{[a;b]}$ and $g \in {\cal U}^{q}\rbr{[a;b]}$ for some $p>1,$ $q>1,$ such that $1/p+1/q>1$ and they have no common points of discontinuity. Then for the constant $E_{p,q}=(p-1)^{1-1/p}p^{-1} D_{p,q}\leq 2 \cdot D_{p,q}$  one has
\[ \left\Vert \int_{a}^{\cdot}\left[f\left(s\right)-f\left(a\right)\right]\mathrm{d}g\left(s\right)  \right\Vert_{q-\emph{TV},\left[a;b\right]} 
\leq E_{p,q}\left\Vert f\right\Vert _{p-\emph{TV},\left[a;b\right]}\left\Vert g\right\Vert _{q-\emph{TV},\left[a;b\right]}.
\]
\end{coro}

\section{Integral equations driven by functions from ${\cal U}^p\rbr{[a;b]}$}
Let $p \in (1;2).$ The preceding section provides us with tools to solve integral equations of the following form
\begin{equation} \label{eq:integral}
y(t)=y_0 + \int_a^t F(y(s)) \dd x(s),
\end{equation}
where $x$ is a continuous function from the space ${\cal U}^p\rbr{[a;b]}$  and $F:\R \ra \R$ is $\alpha$-Lipschitz. 
For our purposes it will be enough to work with the following definition
of locally or globally $\alpha$-Lipschitz function when $\alpha\in\left(0;2\right]$.
For $x=\left(x_{1},\ldots,x_{n}\right)\in\R^{n}$ we denote $\left\Vert x\right\Vert =\max_{i=1,\ldots,n}\left|x_{i}\right|.$
\begin{defi}
Let $F:\R^{n}\ra\R$ and $\alpha\in\left(0;1\right].$
For any $R>0$ we define its \emph{local $\alpha$-Lipschitz parameter
$K_{F}^{\alpha}\left(R\right)$ as} 
\[
K_{F}^{\left(\alpha\right)}\left(R\right):=\sup\left\{ \frac{\left|F\left(y\right)-F\left(x\right)\right|}{\left\Vert y-x\right\Vert ^{\alpha}}:x,y\in\R^{n},x\neq y,\left\Vert x\right\Vert \leq R,\left\Vert y\right\Vert \leq R\right\} 
\]
and its \emph{global $\alpha$-Lipschitz parameter $K_{F}^{\left(\alpha\right)}$
as }$K_{F}^{\left(\alpha\right)}:=\lim_{R\ra+\ns}K_{F}\left(R\right)<+\ns.$
The function $F$ will be called \emph{locally $\alpha$-Lipschitz}
if for every $R>0,$ $K_{F}^{\left(\alpha\right)}\left(R\right)<+\ns$
and it will be called \emph{globally $\alpha$-Lipschitz}
if $K_{F}^{\left(\alpha\right)}<+\ns.$ 
\end{defi}
In the case, when there is no ambiguity what is the value of the parameter
$\alpha$ and what is the function $F,$ we will write $K_{F}\left(R\right),$
$K_{F}$ or even $K\left(R\right)$ or $K.$

First we will consider the case $p-1 < \alpha < 1.$ In this case we have the existence but no uniqueness result. We will obtain a stronger result than similar results \cite[Lemma, p. 459]{Lyons:1994} or \cite[Theorem 1.20]{LyonsCaruana:2007}. Namely, we will prove that there exists a solution to (\ref{eq:integral}) {\bf which is an element of the space} ${\cal U}^p\rbr{[a;b]},$ not only an element of the space ${\cal U}^{q}\rbr{[a;b]}$ for arbitrary chosen $q>p.$ This will be possible with the use of Theorem \ref{thm_integral}.
\begin{prop} \label{prop:int_equation}
Let $p\in(1;2),$ $y_0\in \R,$ $x$ be a continuous function from the space ${\cal U}^p\rbr{[a;b]}$ and $F:\R\ra \R$ be globally $\alpha$-Lipschitz where $p-1 < \alpha < 1.$ Equation (\ref{eq:integral}) admits a solution $y,$ which is an element of ${\cal U}^p\rbr{[a;b]}.$ Moreover, $\TVfullnorm{y}{p}\leq R,$ where $R>0$ satisfies the equality 
\[
R = \rbr{E_{p/\alpha,p} +1}K^{(\alpha)}_F\TVnorm {x}{p} R^{\alpha} + \left| y_0 \right| +\left| F(0)\right|\TVnorm {x}{p},
\]
with $E_{p/\alpha,p}$ being the same as in Corollary \ref{coro_Love_Young_TV}.
\end{prop}
\begin{rem}
Proposition \ref{prop:int_equation} remains true when one exchanges in its formulation the space ${\cal U}^p\rbr{[a;b]}$ for the space space ${\cal V}^p\rbr{[a;b]}.$ The proof goes exactly in the same way, but instead of Theorem \ref{thm_integral} one uses Remark \ref{rem_int_V^p}. 
\end{rem}
First we will prove an auxiliary lemma (analog of \cite[Lemma 1.18]{LyonsCaruana:2007}).
\begin{lema} \label{lema:norm}
Let $p \geq 1,$ $\alpha \in (0;1],$ $f \in {\cal U}^p\rbr{[a;b]}$ and $F:\R \ra \R$ be globally $\alpha$-Lipschitz. Then $F(f(\cdot)) \in {\cal U}^{p/\alpha}\rbr{[a;b]}$ and $\left\Vert F(f(\cdot)) \right\Vert_{p/\alpha-\emph{TV},\left[a;b\right]} \leq K_F^{(\alpha)} \left\Vert f \right\Vert_{p-\emph{TV},\left[a;b\right]}^{\alpha}.$
\end{lema}
\begin{dwd}
For $\alpha \leq 1,$ $x \in \R,$ $K\geq 0$ and $\delta >0$ we have the following elementary inequality 
\begin{equation} \label{ineq_elementary}
\rbr{K|x|^{\alpha}-\delta}_+ \leq K^{1/\alpha} \delta ^{1-1/\alpha} \rbr{|x|-\rbr{\delta/K}^{1/\alpha}}_+.
\end{equation}
Using this inequality and the estimate $\left|F(u) - F(v)\right| \leq K \left|u - v\right|^{\alpha},$ for $K=K_F^{(\alpha)},$ $s,t\in[a;b],$ $\delta > 0$ and $\eta = (\delta/K)^{1/\alpha}$  we calculate 
\begin{eqnarray*}
\delta^{p/\alpha-1}\left(\left|F\left(f\left(t\right)\right)-F\left(f\left(s\right)\right)\right|-\delta\right)_{+} & \leq & \delta^{p/\alpha-1}\left(K\left|f\left(t\right)-f\left(s\right)\right|^{\alpha}-\delta\right)_{+}\\
 & \leq & K^{1/\alpha}\delta^{p/\alpha-1}\delta^{1-1/\alpha}\left(\left|f\left(t\right)-f\left(s\right)\right|-\left(\delta/K\right)^{1/\alpha}\right)_{+}\\
 & = & K^{p/\alpha}\left[\left(\delta/K\right)^{1/\alpha}\right]^{p-1}\left(\left|f\left(t\right)-f\left(s\right)\right|-\left(\delta/K\right)^{1/\alpha}\right)_{+}\\
 & = & K^{p/\alpha}\eta^{p-1}\left(\left|f\left(t\right)-f\left(s\right)\right|-\eta\right)_{+}.
\end{eqnarray*}
From which we get
\[
\sup_{\delta>0}\delta^{p/\alpha-1}\TTV{F\left(f\left(\cdot\right)\right)}{[a;b]}{\delta}\leq K^{p/\alpha}\sup_{\eta>0}\eta^{p-1}\TTV{f}{[a;b]}{\eta},
\]
thus $F(f(\cdot)) \in {\cal U}^{p/\alpha}\rbr{[a;b]}$ and
\[
\left\Vert F(f(\cdot)) \right\Vert_{p/\alpha-\text{TV},\left[a;b\right]} \leq K \left\Vert f \right\Vert_{p-\text{TV},\left[a;b\right]}^{\alpha}.
\]
\end{dwd}
Now we proceed to the proof of Proposition \ref{prop:int_equation}. We will proceed in a standard way, but with the more accurate estimate of Theorem \ref{thm_integral} we will be able to obtain the finiteness of $\TVfullnorm {\cdot}{p}$ norm of the solution.
\begin{dwd}
Let $f \in {\cal U}^{p}\rbr{[a;b]}.$ By Lemma \ref{lema:norm}, $F(f(\cdot)) \in {\cal U}^{p/\alpha}\rbr{[a;b]}$ and since $\alpha/p+1/p>1,$ we may apply Corollary \ref{coro_Love_Young_TV} and define the operator $T: {\cal U}^{p}\rbr{[a;b]} \ra {\cal U}^{p}\rbr{[a;b]},$ 
\[
Tf := y_0 + \int_a^{\cdot} F(f(t))\dd x(t).
\] 
Denote $K=K^{(\alpha)}_F.$ We estimate 
\begin{align}
& \TVfullnorm{Tf}{p} = \TVfullnorm {y_0 + \int_a^{\cdot} F(f(t))\dd x(t)}{p} \nonumber \\ 
& \leq \left| y_0 \right| + \TVnorm {\int_a^{\cdot} \sbr{ F(f(t)) - F(f(a)) }\dd x(t)}{p} + \left| F(f(a))\right|\TVnorm {x}{p} \nonumber \\
& \leq \left| y_0 \right| + \rbr{ E_{p/\alpha,p} \TVnorm {F(f(\cdot))}{p/\alpha} + \left| F(f(a))\right| } \TVnorm {x}{p} \nonumber \\ 
& \leq \left| y_0 \right| + \rbr{ E_{p/\alpha,p} K \TVnorm {f}{p}^{\alpha} + \left| F(f(a))\right| }\TVnorm {x}{p}. \label{ineq:TVnorm_est}
\end{align} 
We naturally have $\TVnorm {f}{p}^{\alpha} \leq \TVfullnorm {f}{p}^{\alpha}$ and, by the Lipschitz property, 
\[
\left| F(f(a))\right| \leq K\left| f(a)\right|^{\alpha} + \left| F(0)\right| \leq K\TVfullnorm {f}{p}^{\alpha} + \left| F(0)\right|.
\]
Denoting 
\begin{equation} \label{def:A}
A = \rbr{E_{p/\alpha,p}+1}K\TVnorm {x}{p} \text{ and } B= \left| y_0 \right| +\left| F(0)\right|\TVnorm {x}{p},
\end{equation}
from (\ref{ineq:TVnorm_est}) we get 
\begin{equation} \label{ineq:stabil}
\TVfullnorm {Tf}{p} \leq A \TVfullnorm {f}{p}^{\alpha} + B.
\end{equation}

For $\alpha <1$ let $R$ be the least positive solution of the inequality $R \geq  A \cdot R^\alpha + B$ (i.e. $R = A \cdot R^\alpha + B$). From (\ref{ineq:stabil}) we have that the operator $T$ maps the closed ball ${\cal B}(R) = \cbr{f\in {\cal U}^{p}\rbr{[a;b]} : \TVfullnorm{f}{p} \leq R}$ to itself. 

Now, for $f,g \in {\cal U}^{p}\rbr{[a;b]}$ we are going to investigate the difference $Tf-Tg.$ Using Theorem \ref{thm_integral}, Lemma \ref{lema:norm} and the Lipschitz property we estimate
\begin{align}
& \TVfullnorm{Tf - Tg}{p} = \TVnorm {\int_a^{\cdot} \sbr{ F(f(t)) -F(g(t))}\dd x(t)}{p}  \nonumber \\ 
& \leq \TVnorm {\int_a^{\cdot} \sbr{ F(f(t)) -F(g(t)) - \cbr{F(f(a)) -F(g(a))}}\dd x(t)}{p} \nonumber \\ & \quad +  \left| F(f(a)) -F(g(a)) \right| \TVnorm {x}{p}  \nonumber \\
& \leq D_{p/\alpha,p} \TVnorm {F(f(\cdot))-F(g(\cdot))}{p/\alpha}^{(p-1)/\alpha} \Oscnorm {F(f(\cdot))-F(g(\cdot))}^{(\alpha+1-p)/\alpha} \TVnorm {x}{p} \nonumber \\ & \quad +  \left| F(f(a)) -F(g(a)) \right| \TVnorm {x}{p} \nonumber \\
&  \leq  D_{p/\alpha,p} K \rbr{\TVnorm {f}{p}^{\alpha}+\TVnorm {g}{p}^{\alpha}}^{(p-1)/\alpha} \Oscnorm {f-g}^{\alpha+1-p} \TVnorm {x}{p} \nonumber \\
& \quad +   K \left|f(a) -g(a) \right|^{\alpha}  \TVnorm {x}{p}. \label{lasttt}
\end{align} 
From (\ref{lasttt}) we see that $T$ is continuous. Moreover, 
from the first inequality in Remark \ref{rem_Young1} and the continuity of $x$ we get that functions belonging to the image $T({\cal B}(R))$ are equicontinuous. Let ${\cal U}$ be the closure of the convex hull of $T({\cal B}(R))$ (in the topology induced by the norm $\TVfullnorm{\cdot}{p}$). It is easy to see that functions belonging to ${\cal U}$ are also equicontinuous. Moreover, ${\cal U} \subset {\cal B}(R)$ and $T({\cal U}) \subset {\cal U}.$ Now, let ${\cal V} = T({\cal U}).$ From the equicontinuity of ${\cal U},$ Arzela-Ascoli Theorem and (\ref{lasttt}) we see that the set ${\cal V}$ is compact in the topology induced by the norm $\TVfullnorm{\cdot}{p}.$  Thus, by the fixed-point Theorem of Schauder, we get that there exists a point $y \in {\cal U}$ such that $Ty=y.$
\end{dwd}

Now we will consider the case $\alpha=1.$ From Remark \ref{rem:no_splitting} it follows that we can not apply frequently used technique of truncating the interval $[a;b]$ to a shorter interval $[a;c],$ $c\in (a;b)$ such that the operator $T$ defined in the proof of Proposition \ref{prop:int_equation} is contractant on ${\cal U}^p\rbr{[a;c]}.$  
To deal with this case let us introduce the following definition.
\begin{defi}\label{interval_split}
Let $p>1$ and $x \in {\cal U}^p\rbr{[a;b]}.$ The function $x$ has {\em splitting property in ${\cal U}^p\rbr{[a;b]}$} if for any $\varepsilon>0$ there exist $\delta >0$ such that $\TVnormAthm x{p}{[c;d]}<\varepsilon$ whenever $0\leq d-c \leq \delta$ and  $[c;d] \subset [a;b].$ 
\end{defi}
Notice that if $x \in {\cal V}^p\rbr{[a;b]}$ then $x$ has the splitting property in ${\cal U}^p\rbr{[a;b]}.$
\begin{fact}
Let $p\in(1;2),$ $y_0\in \R,$ $x$ be a continuous function from the space ${\cal U}^p\rbr{[a;b]}$ and $F:\R\ra \R$ be globally $1$-Lipschitz. Equation (\ref{eq:integral}) admits a solution $y,$ which is an element of ${\cal V}^q\rbr{[a;b]}$ for any $q \in (p;\alpha +1).$ If the function $x$ has the splitting property in ${\cal U}^p\rbr{[a;b]}$ then equation (\ref{eq:integral}) admits a solution $y,$ which is an element of ${\cal U}^p\rbr{[a;b]}.$ If $x\in {\cal V}^p\rbr{[a;b]}$ then equation (\ref{eq:integral}) admits a solution $y,$ which is an element of ${\cal V}^p\rbr{[a;b]}.$
\end{fact}
\begin{dwd} Choose $q,q'\in (p;\alpha +1),$ $q'<q.$  By Proposition \ref{prop:tv_p_norm_var_q_norm} we have $x \in {\cal V}^{q'}\rbr{[a;b]}$ and to prove the existence of the solution $y\in {\cal V}^{q}\rbr{[a;b]}$ we may proceed exactly in the same way as in the proof of \cite[Theorem 1.20]{LyonsCaruana:2007}.

To prove the assertion for $x$ which has the splitting property in ${\cal U}^p\rbr{[a;b]},$ we may proceed in a similar way as in the proof of Proposition \ref{prop:int_equation}. The only thing we need is to assure that the inequality $R\geq A\cdot R + B,$ where $A$ and $B$ are defined in display (\ref{def:A}), holds for sufficiently large $R.$ This may be achieved by splitting the interval $[a;b]$ into small intervals, such that $A < 1$ on every of these intervals, and then solve the equation (\ref{eq:integral}) on every of these intervals with the initial condition being equal the terminal value of the solution on the preceding interval.  

Finally, to prove the assertion for $x\in {\cal V}^{p}\rbr{[a;b]}$ we proceed in a similar way as for $x$ which has the splitting property in ${\cal U}^p\rbr{[a;b]},$ but instead of Theorem \ref{thm_integral} one uses Remark \ref{rem_int_V^p}. 
\end{dwd}

Now we will proceed to the case when $F$ is $(1+\alpha)$-Lipschitz and $\alpha +1 \in [p;2].$ In this case we have the uniqueness result.
We will need the definition of $\left(1+\alpha\right)$-Lipschitz
function with $\alpha\in\left(0;1\right].$  
\begin{defi}
Let $F:\R^{n}\ra\R$ and $\alpha\in\left(0;1\right].$
The function $F$ will be called \emph{locally $\left(1+\alpha\right)$-Lipschitz}
if there exist all partial derivatives $\frac{\partial f}{\partial x_{i}}\left(x\right),$
$i=1,\ldots,n,$ at every point $x\in\R,$ and they are locally $\alpha$-Lipschitz.
Similarly, $F$ will be called \emph{globally $\left(1+\alpha\right)$-Lipschitz}
if it is globally $1$-Lipschitz and there exist all partial derivatives $\frac{\partial f}{\partial x_{i}}\left(x\right),$
$i=1,\ldots,n,$ at every point $x\in\R,$ which are globally $\alpha$-Lipschitz.
\end{defi}

For $\alpha\in\left(0;1\right]$ every locally $(1+\alpha)$-Lipschitz
function $F:\R^{n}\ra\R$ is also locally $1$-Lipschitz
function. Moreover, for any such $\alpha$ and every locally (globally) $(1+\alpha)$-Lipschitz
function there exist 
$G:\R^{n}\times\R^{n}\ra\R,$ which is locally (globally) $\alpha$-Lipschitz
respectively, and for every $x,y\in\R^{n}$ 
\[
F\left(y\right)-F\left(x\right)=G\left(y,x\right)\left(y-x\right)
\]
(see \cite[Proposition 1.26(Division Property)]{LyonsCaruana:2007}).
If $F$ is globally $1$-Lipschitz then $G$ is globally bounded.

We have the following analog of \cite[Theorem 1.28]{LyonsCaruana:2007} for the equations driven by the elements of ${\cal U}^{p}\left(\left[a;b\right]\right)$ 
\begin{prop} \label{thm_int_eq_a>1}
Let $p\in\left(1;2\right),$ $\alpha+1\geq p$ and $y_{0}\in\R.$
If $x\in{\cal U}^{p}\left(\left[a;b\right]\right)$ has the splitting property in ${\cal U}^p\rbr{[a;b]}$ and $F:\R\ra\R$
is locally $\left(\alpha+1\right)$-Lipschitz, then either equation
(\ref{eq:integral}) admits the unique solution $y$ which exists on the whole interval $\left[a;b\right]$
and $y\in{\cal U}^{p}\left(\left[a;b\right]\right)$ or there exists
$d\in\left(a;b\right]$ such that the finite limit $\lim_{t\ra d-}y\left(t\right)$
does not exist but $y\in{\cal U}^{p}\left(\left[a;c\right]\right)$
for every $c\in\left(a;d\right).$ If the function $F$ is globally $\rbr{1+\alpha}$-Lipschitz and globally bounded, then the solution exists on the
whole interval $\left[a;b\right].$ Moreover, denoting this solution
by $I_{F}\left(x,y_{0}\right),$ we get a continuous map $I_{F}:{\cal U}^{p}\left(\left[a;b\right]\right)\times\R\ra{\cal U}^{p}\left(\left[a;b\right]\right).$
\end{prop}
\begin{rem}
If $x$ has no splitting property but is continuous and $\alpha +1 > p$ then for any $q\in (p;\alpha +1)$ we have $x \in {\cal V}^q\rbr{[a;b]}$ and the assertions of Proposition \ref{thm_int_eq_a>1} hold, with the space ${\cal U}^{p}\left(\left[a;b\right]\right)$ replaced by the space ${\cal V}^{q}\left(\left[a;b\right]\right).$
\end{rem}
First we will prove an auxiliary Lemma (analog of \cite[Proposition 1.27]{LyonsCaruana:2007}).
\begin{lema} \label{lema_int_eq_a>1}
Assume that $F:\R\ra\R$ is locally $\left(1+\alpha\right)$-Lipschitz
for some $\alpha\in\left(0;1\right].$ Let $p\geq1.$ If $f,g\in{\cal U}^{p}\left(\left[a;b\right]\right)$
and $\left\Vert f\right\Vert _{\ns,\left[a;b\right]}\leq M,$ $\left\Vert g\right\Vert _{\ns,\left[a;b\right]}\leq M,$
then $F\left(f\left(\cdot\right)\right)-F\left(g\left(\cdot\right)\right)\in{\cal U}^{p/\alpha}\left(\left[a;b\right]\right),$
moreover 
\begin{align*}
& \left\Vert F\left(f\left(\cdot\right)\right)-F\left(g\left(\cdot\right)\right)\right\Vert _{p/\alpha-\emph{TV},\left[a;b\right]} \\ & \leq  2\sup_{x,y \in [-M;M]} \left|G(x,y) \right|^{1-\alpha} \left\Vert f-g\right\Vert _{\emph{osc},\left[a;b\right]}^{1-\alpha}\left\Vert f-g\right\Vert _{p-\emph{TV},\left[a;b\right]}^{\alpha}\\
 &  \quad +4 K_{G}^{\left(\alpha\right)}\left(M\right)\left(\left\Vert f\right\Vert _{p-\emph{TV},\left[a;b\right]}^{\alpha}+\left\Vert g\right\Vert _{p-\emph{TV},\left[a;b\right]}^{\alpha}\right)\left\Vert f-g\right\Vert _{\infty,\left[a;b\right]}.
\end{align*}
\end{lema}
\begin{dwd}
Let $G$ be the quotient function of $F.$ For $\delta>0$
and $a\leq s<t\leq b$ we calculate 
\begin{align}
& \left(\left|F\left(f\left(t\right)\right)-F\left(g\left(t\right)\right)-\left\{ F\left(f\left(s\right)\right)-F\left(g\left(s\right)\right)\right\} \right|-\delta\right)_{+}\nonumber \\
& =\left(\left|G\left(f\left(t\right),g\left(t\right)\right)\left(f\left(t\right)-g\left(t\right)\right)-G\left(f\left(s\right),g\left(s\right)\right)\left(f\left(s\right),g\left(s\right)\right)\right|-\delta\right)_{+}\nonumber \\
& \leq\left(\left|G\left(f\left(t\right),g\left(t\right)\right)\right|\left|f\left(t\right)-g\left(t\right)-\left\{ f\left(s\right)-g\left(s\right)\right\} \right|-\delta/2\right)_{+}\label{eq:pierw_skl}\\
& \quad +\left(\left|G\left(f\left(t\right),g\left(t\right)\right)-G\left(f\left(s\right),g\left(s\right)\right)\right|\left|\left(f\left(s\right)-g\left(s\right)\right)\right|-\delta/2\right)_{+}\label{eq:drug_skl}
\end{align}
Next we estimate (\ref{eq:pierw_skl}). Denoting $G_{M}:=\sup_{x,y \in [-M;M]} \left|G(x,y) \right| $ we have (recall that $\left\Vert f\right\Vert _{\ns,\left[a;b\right]}\leq M,$ $\left\Vert g\right\Vert _{\ns,\left[a;b\right]}\leq M$)
\begin{align}
& \delta^{p/\alpha-1}\left(\left|G\left(f\left(t\right),g\left(t\right)\right)\right|\left|f\left(t\right)-g\left(t\right)-\left\{ f\left(s\right)-g\left(s\right)\right\} \right|-\delta/2\right)_{+}\nonumber \\
& \leq\delta^{p/\alpha-1}\left(G_{M}\left|f\left(t\right)-g\left(t\right)-\left\{ f\left(s\right)-g\left(s\right)\right\} \right|-\delta/2\right)_{+} \nonumber \\
& =2^{-1}\delta^{p/\alpha-p}\delta^{p-1}\left(2G_{M}\left|f\left(t\right)-g\left(t\right)-\left\{ f\left(s\right)-g\left(s\right)\right\} \right|-\delta\right)_{+}\nonumber \\
& \leq 2^{-1}\left(2G_{M}\right)^{p/\alpha-p}\left\Vert f-g
\right\Vert_{\text{osc},\left[a;b\right]}^{p/\alpha-p}\delta^{p-1}\left(2G_{M}\left|f\left(t\right)-g\left(t\right)-\left\{ f\left(s\right)-g\left(s\right)\right\} \right|-\delta\right)_{+},\label{eq:pier_pier_skl}
\end{align}
where the last inequality follows from the fact that for $\delta>2G_{M}\left\Vert f-g\right\Vert _{\text{osc},\left[a;b\right]}$
we have 
\[
\left(2G_{M}\left|f\left(t\right)-g\left(t\right)-\left\{ f\left(s\right)-g\left(s\right)\right\} \right|-\delta\right)_{+}=0.
\]
Denote $K(M) = K_{G}^{\left(\alpha\right)}\left(M\right).$ By the Lipschitz property, \[\left|G\left(f\left(t\right),g\left(t\right)\right)-G\left(f\left(s\right),g\left(s\right)\right)\right|\leq K(M)\left(\left|f\left(t\right)-f\left(s\right)\right|^{\alpha}
+\left|g\left(t\right)-g\left(s\right)\right|^{\alpha}\right).\]
Using this and the inequalities (\ref{elementary_estimate}), (\ref{ineq_elementary}) we further estimate (\ref{eq:drug_skl})
\begin{align}
& \delta^{p/\alpha-1}\left(\left|G\left(f\left(t\right),g\left(t\right)\right)-G\left(f\left(s\right),g\left(s\right)\right)\right|\left|\left(f\left(s\right)-g\left(s\right)\right)\right|-\delta/2\right)_{+}\nonumber \\
& \leq4^{-1}\delta^{p/\alpha-1}\left(4K(M)\left|f\left(t\right)-f\left(s\right)\right|^{\alpha}\left\Vert f-g\right\Vert _{\infty,\left[a;b\right]}-\delta\right)_{+}\nonumber \\
& \quad +4^{-1}\delta^{p/\alpha-1}\left(4K(M)\left|g\left(t\right)-g\left(s\right)\right|^{\alpha}\left\Vert f-g\right\Vert _{\infty,\left[a;b\right]}-\delta\right)_{+}\nonumber \\
& \leq4^{-1}\left(4K(M)\left\Vert f-g\right\Vert _{\infty,\left[a;b\right]}\right)^{p/\alpha}\eta^{p-1}\left(\left|f\left(t\right)-f\left(s\right)\right|-\eta\right)_{+}\label{eq:dru_pier}\\
& \quad +4^{-1}\left(4K(M)\left\Vert f-g\right\Vert _{\infty,\left[a;b\right]}\right)^{p/\alpha}\eta^{p-1}\left(\left|g\left(t\right)-g\left(s\right)\right|-\eta\right)_{+},\label{eq:dru_dru}
\end{align}
where $\eta=\delta^{1/\alpha}/\left(4K(M)\left\Vert f-g\right\Vert _{\infty,\left[a;b\right]}\right)^{1/\alpha}.$
By (\ref{eq:pierw_skl}), (\ref{eq:drug_skl}) and (\ref{eq:pier_pier_skl})-(\ref{eq:dru_dru}),
we get
\begin{align*}
& \sup_{\delta>0}\delta^{p/\alpha-1}\TTV{F\left(f\left(\cdot\right)\right)-F\left(g\left(\cdot\right)\right)}{\left[a;b\right]}{\delta} \\ & \leq  2^{-1}\left(2G_{M}\right)^{p/\alpha-p}\left\Vert f-g\right\Vert _{\text{osc},\left[a;b\right]}^{p/\alpha-p}\sup_{\delta>0}\delta^{p-1}\TTV{2G_{a,b}\left(f-g\right)}{\left[a;b\right]}{\delta}\\
 &  \quad +4^{-1}\left(4K(M)\left\Vert f-g\right\Vert _{\infty,\left[a;b\right]}\right)^{p/\alpha}\sup_{\eta>0}\eta^{p-1}\TTV f{\left[a;b\right]}{\eta}\\
 &  \quad +4^{-1}\left(4K(M)\left\Vert f-g\right\Vert _{\infty,\left[a;b\right]}\right)^{p/\alpha}\sup_{\eta>0}\eta^{p-1}\TTV g{\left[a;b\right]}{\eta}.
\end{align*}
Thus 
\begin{align*}
& \left\Vert F\left(f\left(\cdot\right)\right)-F\left(g\left(\cdot\right)\right)\right\Vert _{p/\alpha-TV,\left[a;b\right]}^{p/\alpha} \\ 
& \leq 2^{p/\alpha-1}G_{M}^{p/\alpha-p}\left\Vert f-g\right\Vert _{osc,\left[a;b\right]}^{p/\alpha-p}\left\Vert f-g\right\Vert _{p-TV,\left[a;b\right]}^{p}\\
&  \quad +4^{p/\alpha-1}K(M)^{p/\alpha}\left\Vert f-g\right\Vert _{\infty,\left[a;b\right]}^{p/\alpha}\left(\left\Vert f\right\Vert _{p-TV,\left[a;b\right]}^{p}+\left\Vert g\right\Vert _{p-TV,\left[a;b\right]}^{p}\right).
\end{align*}
\end{dwd}
Now we proceed to the proof of Proposition \ref{thm_int_eq_a>1}.
\begin{dwd}
To prove that the unique solution exists on some interval $\left[a;c\right],$
$c\in\left(a;b\right],$ we may apply the splitting property of $x$ to obtain the contraction property of Picard iterations on $\left[a;c\right]$
 and proceed exactly in the same way as 
in the proof of \cite[Theorem 1.28]{LyonsCaruana:2007}. We have all the ingredients
needed for the proof when the driving function $x$ belongs to the
space ${\cal U}^{p}\left(\left[a;b\right]\right).$ Namely, instead
of \cite[Remark 1.17]{LyonsCaruana:2007} one uses our Theorem \ref{thm_integral}
and instead of \cite[Proposition 1.27]{LyonsCaruana:2007} one uses Lemma \ref{lema_int_eq_a>1}.
We may continue this procedure on a new interval, with the starting
value equal to the terminal value of the solution on the preceding
interval, i.e. $y\left(c\right).$ We may continue in this way and
obtain that the finite limit $\lim_{t\ra d-}y\left(t\right)$ does
not exist as $t\ra d-$ for some $d\in\left(a;b\right]$ or it may
happen that the endpoints of consecutive intervals on which we construct
our solution tend to some $d_{0}\in\left(a;b\right).$ If it happens
that $\lim_{t\ra d_{0}-}y\left(t\right)$ exists, then we may again
continue the previously applied procedure with the starting value
$\lim_{t\ra d_{0}-}y\left(t\right).$ Thus, the only possibility is
that either equation (\ref{eq:integral}) admits a unique solution $y$ on the whole
interval $\left[a;b\right]$ and $y\in{\cal U}^{p}\left(\left[a;b\right]\right)$
or there exists $d\in\left(a;b\right]$ such that the finite limit
$\lim_{t\ra d-}y\left(t\right)$ does not exist.

To prove the second assertion it is enough to notice that if the functions
$F$ and $G$ are globally bounded and globally $1-$Lipschitz, and
$\alpha$-Lipschitz respectively, then the conditions imposed on $c\in\left(a;b\right]$
such that the functional $T:{\cal U}^{p}\left(\left[a;c\right]\right)\ra{\cal U}^{p}\left(\left[a;c\right]\right),$
\[
Tf:=y_{0}+\int_{a}^{\cdot}F\left(f\left(t\right)\right)\dd x(t),
\]
is contractant, do not depend on the starting value $y_{0}.$ Indeed,
in such a case, by Corollary \ref{coro_Love_Young_TV} we
have 
\begin{align*}
\left\Vert Tf \right\Vert_{p-\text{TV},\left[a;c\right]} 
& \leq E_{p,p}\left\Vert F(f(\cdot))\right\Vert _{p-\text{TV},\left[a;c\right]}\left\Vert x\right\Vert _{p-\text{TV},\left[a;c\right]} + \left| F(f(a)) \right| \left\Vert x\right\Vert _{p-\text{TV},\left[a;c\right]}  \\
& \leq  \rbr{ E_{p,p} K_{F}^{(1)} \left\Vert f\right\Vert _{p-\text{TV},\left[a;c\right]} + \left\Vert F\right\Vert_{\infty}  }\left\Vert x\right\Vert _{p-\text{TV},\left[a;c\right]}.
\end{align*}
Thus, for $c$ small enough for the inequality 
\begin{equation} \label{last_stab}
E_{p,p} K_{F}^{(1)} \left\Vert x\right\Vert _{p-\text{TV},\left[a;c\right]} \leq \frac{1}{2}
\end{equation} 
to hold, the functional $T$ stabilizes the set ${\cal{D}}(R) = \cbr{f\in {\cal U}^p\rbr{[a;b]}: \TVnormA f{p}{[a;c]} \leq R},$ i.e. $T\rbr{{\cal{D}}(R)} \subset {\cal{D}}(R)$ for any $R \geq 2 \left\Vert F\right\Vert_{\infty}  \left\Vert x\right\Vert _{p-\text{TV},\left[a;c\right]}.$ Let us fix any such $R.$

Now notice that if $f(a)=g(a),$ then $\left\Vert f-g\right\Vert _{\ns,\left[a;c\right]} \leq \left\Vert f-g\right\Vert _{\text{osc},\left[a;c\right]},$ in such a case, by Lemma \ref{lema_int_eq_a>1} and (\ref{ineq_TV_osc}) we get that 
\begin{align*}
& \left\Vert F\left(f\left(\cdot\right)\right)-F\left(g\left(\cdot\right)\right)\right\Vert _{p/\alpha-\text{TV},\left[a;c\right]} \\ & \leq  4  \rbr{  \left\Vert G \right\Vert_{\ns}^{1-\alpha} 
 +2 K_{G}^{\left(\alpha\right)}\left(\left\Vert f\right\Vert _{p-\text{TV},\left[a;c\right]}^{\alpha}+\left\Vert g\right\Vert _{p-\text{TV},\left[a;c\right]}^{\alpha}\right)} \left\Vert f-g\right\Vert _{p-\text{TV},\left[a;c\right]}.
\end{align*}
Now, applying Corollary \ref{coro_Love_Young_TV} and the assumption $f(a)=g(a)$ we get 
\begin{align*}
& \TVfullnormA{Tf - Tg}{p}{[a;c]} = \TVnorm {\int_a^{\cdot} \sbr{ F(f(t)) -F(g(t))}\dd x(t)}{p}  \nonumber \\ 
& \leq E_{p/\alpha,p} \TVnormA {F(f(\cdot))-F(g(\cdot))}{p/\alpha}{[a;c]} \TVnormA {x}{p}{[a;c]} \nonumber \\ & \quad +  \left| F(f(a)) -F(g(a)) \right| \TVnormA {x}{p}{[a;c]} \nonumber \\
&  = E_{p/\alpha,p} \TVnormA {F(f(\cdot))-F(g(\cdot))}{p/\alpha}{[a;c]} \TVnormA {x}{p}{[a;c]}.
\end{align*} 
Notice, that indeed $Tf(a)=Tg(a)=y_0.$ 
Thus the functional $T$ is contactant on the set ${\cal{D}}(R)$  from the second iteration, whenever $c$ is small enough for the inequality 
\begin{equation} \label{lastt}
4  E_{p/\alpha,p} \rbr{  \left\Vert G \right\Vert_{\ns, [a;c]} 
 +4 K_{G}^{\left(\alpha\right)} R} \TVnormA {x}{p}{[a;c]} < 1
\end{equation}
to hold. 
Thus for $c$ small enough for the inequalities (\ref{last_stab}) and (\ref{lastt}) to hold, we have one unique solution $y\in {\cal{D}}(R)$ on the interval $[a;c].$ Since such a choice of $c$ does not depend on the initial condition $y_0,$ by the splitting property of $x$ we may extend this solution on the whole interval $[a;b].$ Since we could choose $R$ as large as we pleased, this solution must be the unique solution. 

The rest of the proof goes exactly in the same way as in
the proof of \cite[Theorem 1.28]{LyonsCaruana:2007}. 
\end{dwd}

\textbf{Acknowledgments} I would like to thank Raouf Ghomrasni
for drawing my attention to paper \cite{TronelVladimirov:2000} and Professor Terry Lyons for drawing my attention to paper \cite{Lyons:1994}.

\end{document}